\documentclass[12pt]{article}
\usepackage{mathrsfs}
\usepackage{amssymb} \textwidth 150mm \textheight 220mm \oddsidemargin
15pt \evensidemargin 0pt \topmargin 0cm \headsep 0.3cm

\usepackage{amsmath}
\usepackage{graphicx}
\usepackage{indentfirst}
\usepackage{cite}
\usepackage{float}
\usepackage{CJK}

\begin{document}
\title{\Large\bf{Bifurcation of periodic orbits by perturbing $n$-dimensional piecewise smooth differential systems with two switching planes
\thanks{E-mail address: yangjh@mail.bnu.edu.cn, jihua1113@163.com} }}
\author{{Jihua Yang}\\
{\small \it School of Mathematics and Computer Science, Ningxia Normal University,}\\
 {\small \it  Guyuan, 756000, PR China} }
\date{}
\maketitle \baselineskip=0.9\normalbaselineskip \vspace{-3pt}
\noindent
{\bf Abstract}\, In this paper, the general perturbation problem of  piecewise smooth integrable differential systems with two switching planes is considered. Firstly, when the unperturbed system has a family of periodic orbits, the first order Melnikov vector function is derived which can be used to study the number of periodic orbits bifurcated from the period annuli. Then, by using the obtained Melnikov vector function, we get an upper bound of the number of periodic orbits of a concrete $n$-dimensional piecewise smooth differential system.
\vskip 0.2 true cm
\noindent

 \section{Introduction and main results}
 \setcounter{equation}{0}
\renewcommand\theequation{1.\arabic{equation}}

One of the main problems in the qualitative theory of polynomial differential systems is the study of their limit cycles. Two main questions arise in this setting in dimension two: the study of the  the number of limit cycles depending on the degree of the polynomial, and the study of how many limit cycles emerge from the period annuli around a center when one perturbs it inside a given class of differential equations. These problems have been studied intensively in planar smooth or piecewise smooth differential systems, one can see \cite{SY,YYH,TY,HI,LH,GYC,HS,H,X,X18,WHC,XH,X16,LL,YZ17,YZ18,LHR,CLYZ}. Our main aim  is to bring this study to higher dimension.

%There have been an ocean of excellent works of limit cycle bifurcations for planar smooth or piecewise smooth differential systems, one can see $\heartsuit\heartsuit\heartsuit$.
As far as we know, high-dimensional ($n>2$) smooth or piecewise smooth differential systems with a switching plane have not been studied much. Similar to the 2-dimensional case, there are two main methods which can be used to study the number of periodic orbits for high-dimensional differential systems. One is to use the Melnikov function established in \cite{TH,HSB}. The other is to use the averaging method developed in \cite{BL,LMN,H17,HSB}. In \cite{TH}, the authors established a formula for the first order Melnikov vector function of high-dimensional piecewise smooth differential systems which plays a crucial role in the study of the number of periodic orbits. Recently, a new development to high-dimensional case on the upper bound of periodic orbits was given in \cite{HSB}. For more results, see \cite{LTZ,CCL,XHZ,HJG,X,LMZ,CFF} and the references therein. However, because of the complexity of the calculation of Melnikov function of high-dimensional systems, there is no literature to estimate the upper bound of the number of periodic orbits of high-dimensional system by using Melnikov function.

In the present paper, we first give the first order Melnikov vector function of high-dimensional ($n>2$) piecewise smooth differential systems with two switching planes which can be used to determine the number of periodic orbits bifurcated from period annuli. Then, by using the obtained Melnikov vector function, we study  the number of periodic orbits of an $n$-dimensional perturbed differential system.

Consider an $n$-dimensional ($n\geq2$) piecewise smooth differential system
\begin{eqnarray}
\dot{\mathbf{x}}=f^1(\mathbf{x})+\varepsilon g^1(\mathbf{x}),\ x_1\geq0,x_2\geq0,
\end{eqnarray}
\begin{eqnarray}
\dot{\mathbf{x}}=f^2(\mathbf{x})+\varepsilon g^2(\mathbf{x}),\ x_1>0,x_2<0,
\end{eqnarray}
\begin{eqnarray}
\dot{\mathbf{x}}=f^3(\mathbf{x})+\varepsilon g^3(\mathbf{x}),\ x_1<0,x_2<0,
\end{eqnarray}
\begin{eqnarray}
\dot{\mathbf{x}}=f^4(\mathbf{x})+\varepsilon g^4(\mathbf{x}),\ x_1<0,x_2>0,
\end{eqnarray}
or
\begin{eqnarray}
\dot{\mathbf{x}}=f(\mathbf{x})+\varepsilon g(\mathbf{x}),
\end{eqnarray}
where $\mathbf{x}=(x_1,x_2,\cdots,x_n)^T$, $0\leq\varepsilon\ll1$, $f^k(\mathbf{x})$ and $g^k(\mathbf{x})$ are $C^\infty$ vector functions for $k=1,2,3,4$ and
\begin{eqnarray*}
f(\mathbf{x})=\begin{cases}
f^1(\mathbf{x}),\ x_1\geq0,x_2\geq0,\\
f^2(\mathbf{x}),\ x_1>0,x_2<0,\\
f^3(\mathbf{x}),\ x_1<0,x_2<0,\\
f^4(\mathbf{x}),\ x_1<0,x_2>0,\\
\end{cases}
g(\mathbf{x})=\begin{cases}
g^1(\mathbf{x}),\ x_1\geq0,x_2\geq0,\\
g^2(\mathbf{x}),\ x_1>0,x_2<0,\\
g^3(\mathbf{x}),\ x_1<0,x_2<0,\\
g^4(\mathbf{x}),\ x_1<0,x_2>0.\end{cases}\end{eqnarray*}

In order that system (1.5) has a family of periodic orbits near the origin for $\varepsilon=0$, we make the following assumptions:

\vskip 0.2 true cm
\hangafter 1
\hangindent 2.7em
\noindent
{\bf (A1)} Suppose that $U\subset \mathbb{R}^n$ is an open set with $U\cup \{x_1=0,x_2=0\}\neq\emptyset$. System $(1.k)|_{\varepsilon=0}$ has $n-1$ different $C^\infty$ first integrals $H_i^k(\mathbf{x})$, $i=1,2,\cdots,n-1 $, such that for each $\mathbf{x}\in U^k$, the gradient
$$DH^k_1(\mathbf{x}),\ DH^k_2(\mathbf{x}),\ \cdots,\ DH^k_{n-1}(\mathbf{x}),\ k=1,2,3,4 $$
are linearly independent, where
$$\begin{aligned}U^1=\{\mathbf{x}\in U|x_1\ge0,x_2\geq0\},\ U^2=\{\mathbf{x}\in U|x_1>0,x_2<0\},\\
U^3=\{\mathbf{x}\in U|x_1<0,x_2<0\},\ U^4=\{\mathbf{x}\in U|x_1<0,x_2>0\}.\end{aligned}$$

\hangafter 1
\hangindent 2.8em
\noindent
{\bf (A2)} Let $\mathbf{H}^k(\mathbf{x})=(H^k_1(\mathbf{x}), H^k_2(\mathbf{x}), \cdots, H^k_{n-1}(\mathbf{x}))^T$. There exists an open set $G\subset\mathbb{R}^{n-1}$ such that for each $h=(h_1,h_2,\cdots,h_{n-1})^T\in G$, the curves
$$\begin{aligned}&L_h^1=\{\mathbf{x}\in U^1|\mathbf{H}^1(\mathbf{x})=h\},\qquad\qquad L_h^2=\{\mathbf{x}\in U^2|\mathbf{H}^2(\mathbf{x})=\mathbf{H}^2(B(h))\},\\
&L_h^3=\{\mathbf{x}\in U^3|\mathbf{H}^3(\mathbf{x})=\mathbf{H}^3(C(h))\},\ L_h^4=\{\mathbf{x}\in U^4|\mathbf{H}^4(\mathbf{x})=\mathbf{H}^4(D(h))\}\end{aligned}$$
 contain no critical point of $(1.5)|_{\varepsilon=0}$ and have four different end points $A(h)$, $B(h)$, $C(h)$ and $E(h)$ in $U$ satisfying
 $$\begin{aligned}
 &A(h)=\big(0,a_2(h),\cdot,a_n(h)\big)^T,\  B(h)=\big(b_1(h),0,b_3(h),\cdot,b_n(h)\big)^T,\\
 &C(h)=\big(0,c_2(h),\cdot,c_n(h)\big)^T,\  E(h)=\big(e_1(h),0,e_3(h),\cdot,e_n(h)\big)^T.
 \end{aligned}$$
The system (1.1)$|_{\varepsilon=0}$ has an orbital arc $L^1_h$ starting from $A(h)$ and ending at $B(h)$, the system (1.2)$|_{\varepsilon=0}$ has an orbital arc $L^2_h$ starting from $B(h)$ and ending at $C(h)$, the system (1.3)$|_{\varepsilon=0}$ has an orbital arc $L^3_h$ starting from $C(h)$ and ending at $E(h)$, and the system (1.4)$|_{\varepsilon=0}$ has an orbital arc $L^4_h$ starting from $E(h)$ and ending at $A(h)$. Thus, $L_h=L^1_h\cup L^2_h\cup L^3_h\cup L^4_h$ is a periodic orbit of (1.5)$|_{\varepsilon=0}$ surrounding the origin for $h\in G$.

%defined by $\mathbf{H}^1(\mathbf{x})=h,h\in G, x>0,y>0$
\hangafter 1
\hangindent 2.7em
\noindent
{\bf (A3)} The curves $L_h^k,h\in G,k=1,2,3,4$ are not tangent to the switching plane $x_1=0$ ($x_2=0$ resp.) at points $A(h)$ and $C(h)$ ($B(h)$ and $E(h)$ resp.). That is to say, for $h\in G$,
$$\det\frac{\partial( H^k_1,H^k_2,\cdots,H^k_{n-1})}{\partial(x_2,x_3,\cdots,x_n)}\ \Big(\textup{resp.}\ \det\frac{\partial (H^k_1,H^k_2,\cdots,H^k_{n-1})}{\partial(x_1,x_3,\cdots,x_n)}\Big)$$
is not equal to zero at points $A(h)$ and $C(h)$ (resp. $B(h)$ and $E(h)$).

%There exist an interval $\Sigma=(\alpha,\beta)$, and four points $A=(0,a(h))$, $B=(b(h),0)$, $C=(0,c(h))$ and $D=(d(h),0)$ such that for all $h\in \Sigma$
%$$H^1(A)=H^1(B)=h,\   H^2(B)=H^2(C),  \ H^3(C)=H^3(D),  \ H^4(D)=H^4(A)$$
%with $a(h)c(h)<0$ and $b(h)d(h)<0$.
%\vskip 0.1 true cm

%\noindent
%{\bf Assumption (II).} The system (1.5)$|_{\varepsilon=0}$ has an orbital arc $L^1_h$ starting from $A$ and ending at $B$ defined by $H^1(x,y)=h,h\in\Sigma, x>0,y>0$; the system (1.6)$|_{\varepsilon=0}$ has an orbital arc $L^2_h$ starting from $B$ and ending at $C$ defined by $H^2(x,y)=H^2(B),\ x>0,y<0$; the system (1.7)$|_{\varepsilon=0}$ has an orbital arc $L^3_h$ starting from $C$ and ending at $D$ defined by $H^3(x,y)=H^3(C),\ x<0,y<0$, and the system (1.8)$|_{\varepsilon=0}$ has an orbital arc $L^4_h$ starting from $D$ and ending at $A$ defined by $H^4(x,y)=H^4(D),\ x<0,y>0$. Thus, $L_h=L^1_h\cup L^2_h\cup L^3_h\cup L^4_h$ is a periodic orbit of (1.9)$|_{\varepsilon=0}$ surrounding the origin for $h\in \Sigma$.

 \vskip 0.2 true cm

By assumptions $\bf (A1)-(A3)$, $\{L_h,h\in G\}$ is a family of periodic orbits of system (1.5)$|_{\varepsilon=0}$ and each $L_h$ is piecewise smooth. Without loss of generality, we suppose that $L_h$ has a clockwise orientation, as shown in Fig.\,1.
Suppose that the orbit of system (1.1) starting from $A(h)$. Let $B_\varepsilon(h)=(b_{1\varepsilon}(h),0,b_{3\varepsilon}(h),\cdots,b_{n\varepsilon}(h))$ be its first intersection point with the hyperplane $x_2=0$. Let $C_\varepsilon(h)=(0,c_{1\varepsilon}(h),\cdots,c_{n\varepsilon}(h))$ be the first intersection point of the orbit starting from $B_\varepsilon(h)$ of system (1.2) with the hyperplane $x_1=0$. Let $E_\varepsilon(h)=(e_{1\varepsilon}(h),0,e_{3\varepsilon}(h),\cdots,e_{n\varepsilon}(h))$ be the first intersection point of the orbit starting from $C_\varepsilon(h)$ of system (1.3) with the hyperplane $x_2=0$. Let $A_\varepsilon(h)=(0,a_{2\varepsilon}(h),\cdots,,a_{n\varepsilon}(h))$ be the first intersection point of the orbit starting from $E_\varepsilon(h)$ of system (1.4) with the hyperplane $x_1=0$, see Fig.\,2.

By Lemma 1 in \cite{TH}, one knows that $A_\varepsilon(h)$, $B_\varepsilon(h)$, $C_\varepsilon(h)$ and $E_\varepsilon(h)$ are smooth in $\varepsilon$ with $A_\varepsilon(h)|_{\varepsilon=0}=A(h)$. Then one can define
\begin{eqnarray}
H^1(A_\varepsilon(h))-H^1(A(h))=\varepsilon F(h,\varepsilon).
\end{eqnarray}

Similar to the proof of Lemma 2 in \cite{TH}, one has the following proposition.
 \vskip 0.2 true cm

\noindent
{\bf Proposition 1.1.} {\it For each $h_0\in G$ there exists $\varepsilon_0(h_0)>0$ such that $F(h,\varepsilon)\in C^\infty$ for $0\leq\varepsilon\leq\varepsilon_0$, $h\in G$ with $|h-h_0|<\varepsilon_0$. In particular, $F(h,0)\in C^\infty$ for $h\in G$. Moreover, for a given $h_0\in G$, system (1.5) has a periodic orbit near $L_{h_0}$ if and only if $F(h,\varepsilon)$ has a zero in $h$ near $h_0$ for sufficiently small $\varepsilon>0$.}
\vskip 0.2 true cm

 Similar to the smooth system, we call the function $A\rightarrow A_\varepsilon$ the Poincar\'{e} map of system (1.5). Let $M(h)=F(h,0)$. We call $M(h)$ the first order Melnikov function of system (1.5). Hence, under assumptions ${\bf (A1)-(A3)}$, one can see that an isolated zero of $M(h)$ corresponds to a periodic orbit of (1.5). The formula of $M(h)$ is given in the following theorem.
 \vskip 0.2 true cm

\noindent
{\bf Theorem 1.1.} {\it Under the assumptions ${\bf (A1)-(A3)}$, the first order Melnikov vector function ${M}(h)$ of system (1.5) is
\begin{eqnarray}
\begin{aligned}
{M}(h)=&\int_{\widehat{AB}}D\mathbf{H}^1(\mathbf{x})g^1(\mathbf{x})dt+\overline{D\mathbf{H}^1(A)}\Big[\overline{D\mathbf{H}^4(A)}\Big]^{-1}\\
&\times\underline{D\mathbf{H}^4(E)}\Big[\underline{D\mathbf{H}^3(E)}\Big]^{-1}
\overline{D\mathbf{H}^3(C)}\Big[\overline{D\mathbf{H}^2(C)}\Big]^{-1}
\int_{\widehat{BC}}D\mathbf{H}^2(\mathbf{x})g^2(\mathbf{x})dt\\
&+\overline{D\mathbf{H}^1(A)}\Big[\overline{D\mathbf{H}^4(A)}\Big]^{-1}\underline{D\mathbf{H}^4(E)}\Big[\underline{D\mathbf{H}^3(E)}\Big]^{-1}
\int_{\widehat{CE}}D\mathbf{H}^3(\mathbf{x})g^3(\mathbf{x})dt\\
&+\overline{D\mathbf{H}^1(A)}\Big[\overline{D\mathbf{H}^4(A)}\Big]^{-1}
\int_{\widehat{EA}}D\mathbf{H}^4(\mathbf{x})g^4(\mathbf{x})dt.
\end{aligned}
\end{eqnarray}
Further, if $M(h_0)=0$ and $\det DM(h_0)\neq0$ for some $h_0\in G$, then for $|\varepsilon|$ small enough there exists a unique periodic orbits near $L_{h_0}$ for system (1.5). }
\vskip 0.2 true cm

\noindent
{\bf Remark 1.1.} (i) Let $\Phi$ be a matrix or column vector. Then $\overline{\Phi}$ denotes the matrix (resp. column vector) $\Phi$ removing the first column (resp. first element) and $\underline{\Phi}$ denotes the matrix (resp. column vector) $\Phi$ removing the second column (resp. second element).
\vskip 0.2 true cm

\noindent
(ii) If in (1.1) $\mathbf{H}^1(\mathbf{x})=\mathbf{H}^2(\mathbf{x})$, $\mathbf{H}^3(\mathbf{x})=\mathbf{H}^4(\mathbf{x})$, $f^1(\mathbf{x})=f^2(\mathbf{x})$, $f^3(\mathbf{x})=f^4(\mathbf{x})$,  $g^1(\mathbf{x})=g^2(\mathbf{x})$ and $g^3(\mathbf{x})=g^4(\mathbf{x})$, then the first order Melnikov function $M(h)$ in (1.6) coincides with Theorem 1 obtained in \cite{TH}.
\vskip 0.2 true cm

Next we apply Theorem 1.1 to the following $n$-dimensional piecewise smooth differential system having the form
\begin{eqnarray}\begin{aligned}
&\begin{cases}
\dot{x}_1=x_2+\varepsilon g^1_1(\mathbf{x}),\\
\dot{x}_2=-x_1+\varepsilon g^1_2(\mathbf{x}),\\
\dot{x}_3=\varepsilon g^1_3(\mathbf{x}),\\
\cdots\\
\dot{x}_n=\varepsilon g^1_n(\mathbf{x}),\\
\end{cases}x_1\geq0,x_2\geq0,\quad
\begin{cases}
\dot{x}_1=x_2+\varepsilon g^2_1(\mathbf{x}),\\
\dot{x}_2=-x_1+\varepsilon g^2_2(\mathbf{x}),\\
\dot{x}_3=\varepsilon g^2_3(\mathbf{x}),\\
\cdots\\
\dot{x}_n=\varepsilon g^2_n(\mathbf{x}),\\
\end{cases}x_1>0,x_2<0,\\
&\begin{cases}
\dot{x}_1=x_2+\varepsilon g^3_1(\mathbf{x}),\\
\dot{x}_2=-x_1+\varepsilon g^3_2(\mathbf{x}),\\
\dot{x}_3=\varepsilon g^3_3(\mathbf{x}),\\
\cdots\\
\dot{x}_n=\varepsilon g^3_n(\mathbf{x}),\\
\end{cases}x_1<0,x_2<0,\quad
\begin{cases}
\dot{x}_1=x_2+\varepsilon g^4_1(\mathbf{x}),\\
\dot{x}_2=-x_1+\varepsilon g^4_2(\mathbf{x}),\\
\dot{x}_3=\varepsilon g^4_3(\mathbf{x}),\\
\cdots\\
\dot{x}_n=\varepsilon g^4_n(\mathbf{x}),\\
\end{cases}x_1<0,x_2>0,\\
\end{aligned}\end{eqnarray}
where $\mathbf{x}=(x_1,x_2,\cdots,x_n)^T$,
$$\begin{aligned}
&g^1_i(\mathbf{x})=\sum\limits_{k_1+k_2+\cdots+k_n=0}^ma^i_{k_1k_2\cdots k_n}x_1^{k_1}x_2^{k_2}\cdots x_n^{k_n},\\
&g^2_i(\mathbf{x})=\sum\limits_{k_1+k_2+\cdots+k_n=0}^mb^i_{k_1k_2\cdots k_n}x_1^{k_1}x_2^{k_2}\cdots x_n^{k_n},\\
&g^3_i(\mathbf{x})=\sum\limits_{k_1+k_2+\cdots+k_n=0}^mc^i_{k_1k_2\cdots k_n}x_1^{k_1}x_2^{k_2}\cdots x_n^{k_n},\\
&g^4_i(\mathbf{x})=\sum\limits_{k_1+k_2+\cdots+k_n=0}^md^i_{k_1k_2\cdots k_n}x_1^{k_1}x_2^{k_2}\cdots x_n^{k_n},\\
&i=1,2,\cdots,n.
\end{aligned}$$
Applying the above first order Melnikov function in (1.7), we obtain an upper bound of the number of periodic orbits which bifurcate from the period annulus around origin of system $(1.8)|_{\varepsilon=0}$.
\vskip 0.2 true cm

\noindent
 {\bf Theorem 1.2.} {\it Consider system (1.8) with $\varepsilon>0$ small enough. Using the first order Melnikov function (1.7),
system (1.8) has at most $m^{n-1}$ periodic orbits for $m\geq2$. If $m=1$, then system (1.8) can have 1 periodic orbit.}
\vskip 0.2 true cm

\noindent
{\bf Remark 1.2.} If the switching plane is $x_1=0$, the authors in \cite{LTZ} obtained that the upper bound of the number of periodic orbits is also $m^{n-1}$ by using the averaging method. So, we conjecture that the number of switching planes does not affect the number of periodic orbits.
\vskip 0.2 true cm

The paper is organized as follows. In Section 2, we will prove Theorem 1.1. The proof of Theorem 1.2 will be given in  Section 3. In Appendix, we give a Lemma which will be used in Section 2. It should be noted that the idea of the proof Theorem 1.1 comes from that of H. Tian and  M. Han \cite{TH}.

\section{Proof of Theorem 1.1}
 \setcounter{equation}{0}
\renewcommand\theequation{2.\arabic{equation}}

It is easy to get that
$$\begin{aligned}
\mathbf{H}^1(A_\varepsilon)-\mathbf{H}^1(A)=&[\mathbf{H}^1(A_\varepsilon)-\mathbf{H}^4(A_\varepsilon)]
+[\mathbf{H}^4(A_\varepsilon)-\mathbf{H}^4(E_\varepsilon)]\\
&+[\mathbf{H}^4(E_\varepsilon)-\mathbf{H}^3(E_\varepsilon)]+
[\mathbf{H}^3(E_\varepsilon)-\mathbf{H}^3(C_\varepsilon)]\\
&+[\mathbf{H}^3(C_\varepsilon)-\mathbf{H}^2(C_\varepsilon)]+[\mathbf{H}^2(C_\varepsilon)-\mathbf{H}^2(B_\varepsilon)]\\
&+[\mathbf{H}^2(B_\varepsilon)-\mathbf{H}^1(B_\varepsilon)]+[\mathbf{H}^1(B_\varepsilon)-\mathbf{H}^1(A)]\\
:=&L_1+L_2+L_3+L_4+L_5+L_6+L_7+L_8,
\end{aligned}$$
which follows directly
\begin{eqnarray}
\begin{aligned}
&D_\varepsilon L_1|_{\varepsilon=0}=\Big[D\mathbf{H}^1(A)-D\mathbf{H}^4(A)\Big]D_\varepsilon A_\varepsilon|_{\varepsilon=0},\\
&D_\varepsilon L_2|_{\varepsilon=0}=D\mathbf{H}^4(A)D_\varepsilon A_\varepsilon|_{\varepsilon=0}
-D\mathbf{H}^4(E)\Big]D_\varepsilon E_\varepsilon|_{\varepsilon=0},\\
&D_\varepsilon L_3|_{\varepsilon=0}=\Big[D\mathbf{H}^4(E)-D\mathbf{H}^3(E)\Big]D_\varepsilon E_\varepsilon|_{\varepsilon=0},\\
&D_\varepsilon L_4|_{\varepsilon=0}=D\mathbf{H}^3(E)D_\varepsilon E_\varepsilon|_{\varepsilon=0}
-D\mathbf{H}^3(C)\Big]D_\varepsilon C_\varepsilon|_{\varepsilon=0},\\
&D_\varepsilon L_5|_{\varepsilon=0}=\Big[D\mathbf{H}^3(C)-D\mathbf{H}^2(C)\Big]D_\varepsilon C_\varepsilon|_{\varepsilon=0},\\
&D_\varepsilon L_6|_{\varepsilon=0}=D\mathbf{H}^2(C)D_\varepsilon C_\varepsilon|_{\varepsilon=0}
-D\mathbf{H}^2(B)\Big]D_\varepsilon B_\varepsilon|_{\varepsilon=0},\\
&D_\varepsilon L_7|_{\varepsilon=0}=\Big[D\mathbf{H}^2(B)-D\mathbf{H}^1(B)\Big]D_\varepsilon B_\varepsilon|_{\varepsilon=0},\\
&D_\varepsilon L_8|_{\varepsilon=0}=D\mathbf{H}^1(B)D_\varepsilon B_\varepsilon|_{\varepsilon=0}.\\
\end{aligned}
\end{eqnarray}
Notice that the first (resp. second) component of both $D_\varepsilon A_\varepsilon$ and $D_\varepsilon C_\varepsilon$ (resp. $D_\varepsilon B_\varepsilon$ and $D_\varepsilon E_\varepsilon$) equal to 0. We have for $i=1,2,3,4$
\begin{eqnarray}
\begin{aligned}
&D\mathbf{H}^i(A)D_\varepsilon A_\varepsilon=\overline{D\mathbf{H}^i(A)}\overline{D_\varepsilon A_\varepsilon},\
D\mathbf{H}^i(C)D_\varepsilon C_\varepsilon=\overline{D\mathbf{H}^i(C)}\overline{D_\varepsilon C_\varepsilon},\\
&D\mathbf{H}^i(B)D_\varepsilon B_\varepsilon=\underline{D\mathbf{H}^i(B)}\underline{D_\varepsilon B_\varepsilon},\
D\mathbf{H}^i(E)D_\varepsilon E_\varepsilon=\underline{D\mathbf{H}^i(E)}\underline{D_\varepsilon E_\varepsilon}.\\
\end{aligned}
\end{eqnarray}
By assumption {\bf (A3)} the square matrices $\overline{D\mathbf{H}^i(A)}$, $\underline{D\mathbf{H}^i(B)}$, $\overline{D\mathbf{H}^i(C)}$ and $\underline{D\mathbf{H}^i(E)}$ are invertible.
From assumption {\bf (A1)}, one knows that $D\mathbf{H}^k(\mathbf{x})f^k(\mathbf{x})=0,k=1,2,3,4$. Hence,
$$\begin{aligned}
L_8=&\mathbf{H}^1(B_\varepsilon)-\mathbf{H}^1(A)=\int_{\widehat{AB_\varepsilon}}d\mathbf{H}^1\\
=&\int_{\widehat{AB_\varepsilon}}D\mathbf{H}^1(\mathbf{x})[f^1(\mathbf{x})+\varepsilon g^1(\mathbf{x})]dt\\
=&\varepsilon\int_{\widehat{AB_\varepsilon}}D\mathbf{H}^1(\mathbf{x})g^1(\mathbf{x})dt\\
=&\varepsilon\int_{\widehat{AB}}D\mathbf{H}^1(\mathbf{x})g^1(\mathbf{x})dt+O(\varepsilon^2),
\end{aligned}$$
which implies
\begin{eqnarray}
D_\varepsilon L_8|_{\varepsilon=0}=\int_{\widehat{AB}}D\mathbf{H}^1(\mathbf{x})g^1(\mathbf{x})dt.
\end{eqnarray}
Similarly, one can get
\begin{eqnarray}
\begin{aligned}
&D_\varepsilon L_6|_{\varepsilon=0}=\int_{\widehat{BC}}D\mathbf{H}^2(\mathbf{x})g^2(\mathbf{x})dt,\\
&D_\varepsilon L_4|_{\varepsilon=0}=\int_{\widehat{CE}}D\mathbf{H}^3(\mathbf{x})g^3(\mathbf{x})dt,\\
&D_\varepsilon L_2|_{\varepsilon=0}=\int_{\widehat{EA}}D\mathbf{H}^4(\mathbf{x})g^4(\mathbf{x})dt.\\
\end{aligned}
\end{eqnarray}
In view of (2.1)-(2.3), one obtains
\begin{eqnarray}
\underline{D_\varepsilon B_\varepsilon}|_{\varepsilon=0}=
\Big[\underline{D\mathbf{H}^1(B)}\Big]^{-1}\int_{\widehat{AB}}D\mathbf{H}^1(\mathbf{x})g^1(\mathbf{x})dt.
\end{eqnarray}
From (2.1), (2.2), (2.4) and (2.5) one has
\begin{eqnarray}
\begin{aligned}
\overline{D_\varepsilon C_\varepsilon}|_{\varepsilon=0}=&
\Big[\overline{D\mathbf{H}^2(C)}\Big]^{-1}\int_{\widehat{BC}}D\mathbf{H}^2(\mathbf{x})g^2(\mathbf{x})dt\\&+
\Big[\overline{D\mathbf{H}^2(C)}\Big]^{-1}\underline{D\mathbf{H}^2(B)}\Big[\underline{D\mathbf{H}^1(B)}\Big]^{-1}
\int_{\widehat{AB}}D\mathbf{H}^1(\mathbf{x})g^1(\mathbf{x})dt,\\
\overline{D_\varepsilon E_\varepsilon}|_{\varepsilon=0}=&
\Big[\underline{D\mathbf{H}^3(E)}\Big]^{-1}\int_{\widehat{CE}}D\mathbf{H}^3(\mathbf{x})g^3(\mathbf{x})dt\\&+
\Big[\underline{D\mathbf{H}^3(E)}\Big]^{-1}\overline{D\mathbf{H}^3(C)}\Big[\overline{D\mathbf{H}^2(C)}\Big]^{-1}
\int_{\widehat{BC}}D\mathbf{H}^2(\mathbf{x})g^2(\mathbf{x})dt\\
&+\Big[\underline{D\mathbf{H}^3(E)}\Big]^{-1}\overline{D\mathbf{H}^3(C)}\Big[\overline{D\mathbf{H}^2(C)}\Big]^{-1}
\underline{D\mathbf{H}^2(B)}\Big[\underline{D\mathbf{H}^1(B)}\Big]^{-1}\\
&\times\int_{\widehat{AB}}D\mathbf{H}^1(\mathbf{x})g^1(\mathbf{x})dt,\\
\overline{D_\varepsilon A_\varepsilon}|_{\varepsilon=0}=&\Big[\overline{D\mathbf{H}^4(A)}\Big]^{-1}\int_{\widehat{EA}}D\mathbf{H}^4(\mathbf{x})g^4(\mathbf{x})dt\\
&+\Big[\overline{D\mathbf{H}^4(A)}\Big]^{-1} \underline{D\mathbf{H}^4(E)} \Big[\underline{D\mathbf{H}^3(E)}\Big]^{-1}\int_{\widehat{CE}}D\mathbf{H}^3(\mathbf{x})g^3(\mathbf{x})dt\\
&+\Big[\overline{D\mathbf{H}^4(A)}\Big]^{-1} \underline{D\mathbf{H}^4(E)}\Big[\underline{D\mathbf{H}^3(E)}\Big]^{-1}\overline{D\mathbf{H}^3(C)}\Big[\overline{D\mathbf{H}^2(C)}\Big]^{-1}\\
&\times\int_{\widehat{BC}}D\mathbf{H}^2(\mathbf{x})g^2(\mathbf{x})dt\\
&+\Big[\overline{D\mathbf{H}^4(A)}\Big]^{-1} \underline{D\mathbf{H}^4(E)}\Big[\underline{D\mathbf{H}^3(E)}\Big]^{-1}\overline{D\mathbf{H}^3(C)}\Big[\overline{D\mathbf{H}^2(C)}\Big]^{-1}\\
&\times\underline{D\mathbf{H}^2(B)}\Big[\underline{D\mathbf{H}^1(B)}\Big]^{-1}\int_{\widehat{AB}}D\mathbf{H}^1(\mathbf{x})g^1(\mathbf{x})dt,\\
\end{aligned}
\end{eqnarray}
By (1.6), one has
$$D_\varepsilon \big[\mathbf{H}^1(A_\varepsilon(h))-\mathbf{H}^1(A(h))\big]=F(h,\varepsilon)+\varepsilon D_\varepsilon F(h,\varepsilon).$$
Hence, it follows from $M(h)=F(h,0)$ that
\begin{eqnarray}
M(h)=\sum\limits_{i=1}^8D_\varepsilon L_i|_{\varepsilon=0}.
\end{eqnarray}
Now, combining with (2.1), (2.3), (2.4) and (2.6), we can get
\begin{eqnarray}
\begin{aligned}
\mathbf{M}(h)=&\overline{D\mathbf{H}^1(A)}\Big[\overline{D\mathbf{H}^4(A)}\Big]^{-1}\underline{D\mathbf{H}^4(E)}\Big[\underline{D\mathbf{H}^3(E)}\Big]^{-1}
\overline{D\mathbf{H}^3(C)}\Big[\overline{D\mathbf{H}^2(C)}\Big]^{-1}\\
&\times\underline{D\mathbf{H}^2(B)}\Big[\underline{D\mathbf{H}^1(B)}\Big]^{-1}\int_{\widehat{AB}}D\mathbf{H}^1(\mathbf{x})g^1(\mathbf{x})dt\\
&+\overline{D\mathbf{H}^1(A)}\Big[\overline{D\mathbf{H}^4(A)}\Big]^{-1}\underline{D\mathbf{H}^4(E)}\Big[\underline{D\mathbf{H}^3(E)}\Big]^{-1}\\
&\times\overline{D\mathbf{H}^3(C)}\Big[\overline{D\mathbf{H}^2(C)}\Big]^{-1}
\int_{\widehat{BC}}D\mathbf{H}^2(\mathbf{x})g^2(\mathbf{x})dt\\
&+\overline{D\mathbf{H}^1(A)}\Big[\overline{D\mathbf{H}^4(A)}\Big]^{-1}\underline{D\mathbf{H}^4(E)}\Big[\underline{D\mathbf{H}^3(E)}\Big]^{-1}
\int_{\widehat{CE}}D\mathbf{H}^3(\mathbf{x})g^3(\mathbf{x})dt\\
&+\overline{D\mathbf{H}^1(A)}\Big[\overline{D\mathbf{H}^4(A)}\Big]^{-1}
\int_{\widehat{EA}}D\mathbf{H}^4(\mathbf{x})g^4(\mathbf{x})dt.
\end{aligned}
\end{eqnarray}
From assumption {\bf (A2)}, we have
$$\begin{aligned}
&\mathbf{H}^1(A(h))=\mathbf{H}^1(B(h))=h,\ \mathbf{H}^2(B(h))=\mathbf{H}^2(C(h)),\\
&\mathbf{H}^3(C(h))=\mathbf{H}^3(E(h)),\ \mathbf{H}^4(E(h))=\mathbf{H}^4(A(h)).
\end{aligned}$$
Differentiating both sides of the above four equalities with respect to $h$ yields
\begin{eqnarray}
\overline{D\mathbf{H}^1(A(h))}\Big[\overline{[DA(h)]^T}\Big]^T=\underline{D\mathbf{H}^1(B(h))}\Big[\underline{[DB(h)]^T}\Big]^T=I,
\end{eqnarray}
and
\begin{eqnarray}\begin{aligned}
&\underline{D\mathbf{H}^2(B(h))}\Big[\underline{[DB(h)]^T}\Big]^T=\overline{D\mathbf{H}^2(C(h))}\Big[\overline{[DC(h)]^T}\Big]^T,\\
&\overline{D\mathbf{H}^3(C(h))}\Big[\overline{[DC(h)]^T}\Big]^T=\underline{D\mathbf{H}^3(E(h))}\Big[\underline{[DE(h)]^T}\Big]^T,\\
&\underline{D\mathbf{H}^4(E(h))}\Big[\underline{[DE(h)]^T}\Big]^T=\overline{D\mathbf{H}^4(A(h))}\Big[\overline{[DA(h)]^T}\Big]^T,
\end{aligned}\end{eqnarray}
where $I$ is an $(n-1)\times (n-1)$ identity matrix. From (2.9), we have
\begin{eqnarray}\Big[\underline{[DB(h)]^T}\Big]^T=\Big[\underline{D\mathbf{H}^1(B(h))}\Big]^{-1}.\end{eqnarray}
Hence, by (2.10) and (2.11), it follows that
\begin{eqnarray}\begin{aligned}
\Big[\overline{[DA(h)]^T}\Big]^T=&
\Big[\overline{D\mathbf{H}^4(A)}\Big]^{-1}\underline{D\mathbf{H}^4(E)}\Big[\underline{D\mathbf{H}^3(E)}\Big]^{-1}\\&\times
\overline{D\mathbf{H}^3(C)}\Big[\overline{D\mathbf{H}^2(C)}\Big]^{-1}\underline{D\mathbf{H}^2(B)}\Big[\underline{D\mathbf{H}^1(B)}\Big]^{-1}.
\end{aligned}\end{eqnarray}
Inserting (2.12) into (2.9) gives
{\small $$\overline{D\mathbf{H}^1(A)}\Big[\overline{D\mathbf{H}^4(A)}\Big]^{-1}\underline{D\mathbf{H}^4(E)}\Big[\underline{D\mathbf{H}^3(E)}\Big]^{-1}
\overline{D\mathbf{H}^3(C)}\Big[\overline{D\mathbf{H}^2(C)}\Big]^{-1}\underline{D\mathbf{H}^2(B)}\Big[\underline{D\mathbf{H}^1(B)}\Big]^{-1}=I.$$}
Therefore, (1.7) follows from (2.8) and the above equality. This ends the proof.\quad $\lozenge$
\vskip 0.2 true cm

Then, by using Theorem 1.1, we calculate the first order Melnikov vector function $M(h,\delta)$ of a class of $n$-dimensional piecewise smooth differential systems having the form
\begin{eqnarray}
\begin{cases}
\dot{x}_1=H_{x_2}(x_1,x_2,\mathbf{y})+\varepsilon P(x_1,x_2,\mathbf{y},\delta),\\
\dot{x}_2=-H_{x_1}(x_1,x_2,\mathbf{y})+\varepsilon Q(x_1,x_2,\mathbf{y},\delta),\\
\dot{\mathbf{y}}=\varepsilon R(x_1,x_2,\mathbf{y},\delta),\\
\end{cases}
\end{eqnarray}
where $\mathbf{y}=(x_3,x_4,\cdot,x_{n-2})^T\in \mathbb{R}^{n-2},n\geq2,\delta\in D\subset \mathbb{R}^m$ is a vector parameter with $D$ compact and
\begin{eqnarray}\begin{aligned}
&H(x_1,x_2,\mathbf{y})=\begin{cases}
H^1,\ x_1\geq0,x_2\geq0,\\  H^2,\ x_1>0,x_2<0,\\
H^3,\ x_1<0,x_2<0,\\  H^4,\ x_1<0,x_2>0,\\
\end{cases}\\
&P(x_1,x_2,\mathbf{y},\delta)=\begin{cases}
g^1_1,\ x_1\geq0,x_2\geq0,\\  g_1^2,\ x_1>0,x_2<0,\\
g_1^3,\ x_1<0,x_2<0,\\  g_1^4,\ x_1<0,x_2>0,\\
\end{cases}\\
&Q(x_1,x_2,\mathbf{y},\delta)=\begin{cases}
g^1_2,\ x_1\geq0,x_2\geq0,\\  g_2^2,\ x_1>0,x_2<0,\\
g_2^3,\ x_1<0,x_2<0,\\  g_2^4,\ x_1<0,x_2>0,\\
\end{cases}\\
&R(x_1,x_2,\mathbf{y},\delta)=\begin{cases}
(g_3^1,g_4^1,\cdots,g_n^1)^T,\ x_1\geq0,x_2\geq0,\\  (g_3^2,g_4^2,\cdots,g_n^2)^T,\ x_1>0,x_2<0,\\
(g_3^3,g_4^3,\cdots,g_n^3)^T,\ x_1<0,x_2<0,\\  (g_3^4,g_4^4,\cdots,g_n^4)^T,\ x_1<0,x_2>0,
\end{cases}\end{aligned}\end{eqnarray}
with $H^k,g^k_i(k=1,2,3,4;i=1,2,\cdots,n)$ $C^\infty$ functions.

It is easy to verify that the unperturbed system of (2.13) has $H,x_3,x_4,\cdot,x_n$ as its $n-1$ first integrals. The first two equations of (2.13) define a planar Hamiltonian system
\begin{eqnarray}
\begin{cases}
\dot{x}_1=H_{x_2}(x_1,x_2,\mathbf{y}),\\
\dot{x}_2=-H_{x_1}(x_1,x_2,\mathbf{y})
\end{cases}
\end{eqnarray}
with Hamiltonian function $H(x_1,x_2,\mathbf{y})$ containing $n-2$ parameters $x_i(i=3,4,\cdots,n)$. Now we make two basic assumptions ${\bf (H1)}$ and ${\bf (H2)}$ for system (2.13) corresponding to ${\bf (A1)}$-${\bf (A3)}$.
\vskip 0.2 true cm

 \hangafter 1
\hangindent 2.6em
\noindent
${\bf (H1)}$ For each $\hat{h}\in G_1\subseteq \mathbb{R}^{n-2}$ with $G_1$ an open set, there is an open interval $\hat{J}$ dependent on $\hat{h}$ such that
$$\begin{aligned}
&\widehat{AB}:\ L_{h_1,\hat{h}}^1=\{(x_1,x_2)|H^1(x_1,x_2,\hat{h})=h_1,\ x\geq0,x_2\geq0\},\\
&\widehat{BC}:\ L_{h_1,\hat{h}}^2=\{(x_1,x_2)|H^2(x_1,x_2,\hat{h})=H^2(B(h),\hat{h}),\ x>0,x_2<0\},\\
&\widehat{CE}:\ L_{h_1,\hat{h}}^3=\{(x_1,x_2)|H^3(x_1,x_2,\hat{h})=H^3(C(h),\hat{h}),\ x<0,x_2<0\},\\
\end{aligned}$$
and
$$\widehat{EA}:\ L_{h_1,\hat{h}}^4=\{(x_1,x_2)|H^4(x_1,x_2,\hat{h})=H^4(D(h),\hat{h}),\ x<0,x_2>0\}$$
 contain no critical points of system (2.15). System (2.15) has an orbital arc $L^1_{h_1,\hat{h}}$ starting from $A(h)$ and ending at $B(h)$, an orbital arc $L^2_{h_1,\hat{h}}$ starting from $B(h)$ and ending at $C(h)$, an orbital arc $L^3_{h_1,\hat{h}}$ starting from $C(h)$ and ending at $E(h)$ and  an orbital arc $L^4_{h_1,\hat{h}}$ starting from $E(h)$ and ending at $A(h)$.  Thus, $L_h=L^1_{h_1,\hat{h}}\cup L^2_{h_1,\hat{h}}\cup L^3_{h_1,\hat{h}}\cup L^4_{h_1,\hat{h}}$ is a periodic orbit of (2.15).

\hangafter 1
\hangindent 2.6em
\noindent
${\bf (H2)}$ The curves $L^k_{h_1,\hat{h}}(k=1,2,3,4)$ are not tangent to $x_2$-axis (resp. $x_1$-axis) at points $A(h)$ and $C(h)$ (resp. $B(h)$ and $E(h)$). In other words, for each $\hat{h}\in G_1$, $h_1\in J_{\hat{h}}$,
$$H^k_y\big(A(h),\hat{h}\big)H^k_y\big(C(h),\hat{h}\big)H^k_x\big(B(h),\hat{h}\big)H^k_x\big(E(h),\hat{h}\big) \neq0,\ k=1,2,3,4.$$

Next, we will give the first order Melnikov vector function of system (2.13) by Theorem 1.1. For the sake of simplicity, we take $H^1=H^2=H^3=H^4=H$ in (2.14). In this case, the coefficient matrices of the curvilinear integrals in (1.7) are identity matrices. That is,
\begin{eqnarray}\begin{aligned}
&\overline{D\mathbf{H}^1(A)}\Big[\overline{D\mathbf{H}^4(A)}\Big]^{-1}\underline{D\mathbf{H}^4(E)}\Big[\underline{D\mathbf{H}^3(E)}\Big]^{-1}
\overline{D\mathbf{H}^3(C)}\Big[\overline{D\mathbf{H}^2(C)}\Big]^{-1}=I,\\
&\overline{D\mathbf{H}^1(A)}\Big[\overline{D\mathbf{H}^4(A)}\Big]^{-1}\underline{D\mathbf{H}^4(E)}\Big[\underline{D\mathbf{H}^3(E)}\Big]^{-1}
=I,\\
&\overline{D\mathbf{H}^1(A)}\Big[\overline{D\mathbf{H}^4(A)}\Big]^{-1}=I.
\end{aligned}\end{eqnarray}
\vskip 0.2 true cm

\noindent
{\bf Theorem 2.1.}\, {\it Assume that {\bf (H1)} and {\bf (H2)} are satisfied. Then the first order Melnikov vector function $M(h,\delta)$ of system (2.13) can be written as
\begin{eqnarray}\begin{aligned}
M(h,\delta)=&\left(
  \begin{array}{c}
          \sum\limits_{k=1}^4\big(M_1^k(h,\delta)+N_k(h,\delta)\big)\\
          \sum\limits_{k=1}^4M_2^k(h,\delta)\\
          \vdots \\
          \sum\limits_{k=1}^4M_{n-1}^k(h,\delta)
 \end{array}
 \right)\\
 =&\big(M_1(h,\delta),M_2(h,\delta),\cdots,M_{n-1}(h,\delta)\big)^T,\ h=(h_1,\hat{h})^T,
\end{aligned}\end{eqnarray}
where
\begin{eqnarray}\begin{aligned}
&M_1^k(h,\delta)=\int_{L^k_{h_1,\hat{h}}}g^k_2(x_1,x_2,\hat{h},\delta)dx_1-g^k_1(x_1,x_2,\hat{h},\delta)dx_2,\\
&N_k=\sum\limits_{i=1}^{n-2}\int_{L^k_{h_1,\hat{h}}}H_{x_{i+2}}(x_1,x_2,\hat{h},\delta)g^k_{i+2}(x_1,x_2,\hat{h},\delta)dt,\\
&M^k_j(h,\delta)=\int_{L^k_{h_1,\hat{h}}}g^k_{j+1}(x_1,x_2,\hat{h},\delta)dt,\ k=1,2,3,4;j=2,\cdots,n-1.
\end{aligned}\end{eqnarray}}
\vskip 0.2 true cm

\noindent
{\bf Proof.} Denote $$\begin{aligned}&\mathbf{H}(x_1,x_2,\mathbf{y})=(H,x_3,x_4,\cdots,x_n)^T,\\ &g^k(x_1,x_2,\mathbf{y})=(g^k_1,g_2^k,\cdots,g^k_n)^T.\end{aligned}$$
Then, by straightforward computation, one has for $k=1,2,3,4$
$$\begin{aligned}
D\mathbf{H}(x_1,x_2,\mathbf{y})g^k(x_1,x_2,\mathbf{y})=&\left(
  \begin{matrix}
          H_{x_1}&H_{x_2}&H_{x_3}&H_{x_4}&\cdots&H_{x_n}\\
          0&0&1&0&\cdots&0\\
          0&0&0&1&\cdots&0\\
          \vdots&\vdots&\vdots&\vdots&\vdots&\vdots \\
           0&0&0&0&\cdots&1\\
 \end{matrix}
 \right)
 \left(
  \begin{matrix}
          g^k_1\\
          g^k_2\\
          g^k_3\\
          \vdots \\
           g^k_n\\
 \end{matrix}
 \right)\\
 =&\left(
  \begin{matrix}
         \sum\limits_{i=1}^nH_{x_i}g^k_i\\
          g^k_3\\
          g^k_4\\
          \vdots \\
           g^k_n\\
 \end{matrix}
 \right).
 \end{aligned}$$
By substituting the above equality into (1.7) and in view of (2.16), we obtain (2.17). This completes the proof. \quad $\lozenge$
\vskip 0.2 true cm

Now, we give a lemma which provides an effective way to calculate $M^k_2(h,\delta),M^k_3(h,\delta),\cdots,M^k_{n-1}(h,\delta)$ in (2.18).
\vskip 0.2 true cm

\noindent
{\bf Lemma 2.1.} {\it Let $$\begin{aligned}&\bar{R}_{j+1}^k(x_1,x_2,\hat{h},\delta)=\int_0^{x_2}g^k_{j+1}(x_1,x_2,\hat{h},\delta)dx_2,\\
&\bar{M}_{j}^k(h,\delta)=\int_{L^k_{h_1,\hat{h}}}\bar{R}_{j+1}^k(x_1,x_2,\hat{h},\delta)dx_1.\\
\end{aligned}$$
Then, we obtain
$$M_j^k(h,\delta)=\frac{\partial \bar{M}_{j}^k(h,\delta)}{\partial h_1},\ k=1,2,3,4;j=2,3,\cdots,n-1.$$}
\vskip 0.2 true cm

\noindent
{\bf Proof.} By Lemma A.1 in Appendix, one gets for smooth functions $p$ and $q$
$$\begin{aligned}
&\frac{\partial \Big(\int_{{L^1_{h_1,\hat{h}}}}q(x_1,x_2,\hat{h},\delta)dx_1-p(x_1,x_2,\hat{h},\delta)dx_2\Big)}{\partial h_1}\\
&=\int_{{L^1_{h_1,\hat{h}}}}(p_x+q_y)dt+q(B(h),\hat{h},\delta)\frac{\partial b(h)}{h_1}+p(A(h),\hat{h},\delta)\frac{\partial a(h)}{\partial h_1}.
\end{aligned}$$
Particularly, taking $q=\bar{R}_{j+1}^1$ and $p=0$, one has
$$\frac{\partial \bar{M}_{j}^1(h,\delta)}{\partial h_1}=\int_{{L^1_{h_1,\hat{h}}}}\big(\bar{R}_{j+1}^1\big)_{x_2}dt+\bar{R}_{j+1}^1(B(h),\hat{h},\delta).$$
Notice that $$\bar{R}_{j+1}^1(B(h),\hat{h},\delta)=\bar{R}_{j+1}^1(b(h),0,\hat{h},\delta)=0,$$
we get the desired result. The others can proved similarly. This completes the proof. \quad $\lozenge$
\vskip 0.2 true cm

\noindent
{\bf Remark 2.1.} If $H(x_1,x_2,\mathbf{y})$ is independent of $\mathbf{y}$, then the first order Melnikov function in Theorem 2.1 can be written as
\begin{eqnarray}\begin{aligned}
M(h,\delta)=&\left(
  \begin{array}{c}
          \sum\limits_{k=1}^4M_1^k(h,\delta)\\
          \sum\limits_{k=1}^4M_2^k(h,\delta)\\
          \vdots \\
          \sum\limits_{k=1}^4M_{n-1}^k(h,\delta)
 \end{array}
 \right).
\end{aligned}\end{eqnarray}

\section{An application to an $n$-dimensional differential system}
 \setcounter{equation}{0}
\renewcommand\theequation{3.\arabic{equation}}

In this section, we estimate the number of periodic orbits of the $n$-dimensional piecewise smooth differential system (1.8).
For $\varepsilon=0$, system (1.8) has $n-1$ first integrals
\begin{eqnarray}
H(x_1,x_2)=\frac{1}{2}(x_1^2+x_2^2),\ x_3,\ x_4,\cdots,x_n.
\end{eqnarray}
Apparently, $(0,0,\hat{h})$ is a linear center in the plane
 $\mathbf{y}=\hat{h}$, where
  \begin{eqnarray}\mathbf{y}=(x_3,x_4,\cdots,x_n)^T,\ \ \hat{h}=(h_3,h_4,\cdots,h_n)^T.\end{eqnarray}
 In the remainder of this section, we use the same notions as in the previous section. Thus, $L_{h_1,\hat{h}}=L^1_{h_1,\hat{h}}\cup L^2_{h_1,\hat{h}}\cup L^3_{h_1,\hat{h}} \cup L^4_{h_1,\hat{h}}$ is a family of periodic orbits of the following differential system
\begin{eqnarray}\begin{aligned}
&\begin{cases}
\dot{x}_1=x_2,\\
\dot{x}_2=-x_1,\\
\end{cases}x_1\geq0,x_2\geq0,\quad
\begin{cases}
\dot{x}_1=x_2,\\
\dot{x}_2=-x_1,\\
\end{cases}x_1>0,x_2<0,\\
&\begin{cases}
\dot{x}_1=x_2,\\
\dot{x}_2=-x_1,\\
\end{cases}x_1<0,x_2<0,\quad
\begin{cases}
\dot{x}_1=x_2,\\
\dot{x}_2=-x_1,\\
\end{cases}x_1<0,x_2>0,\\
\end{aligned}\end{eqnarray}
where
$$\begin{aligned}
&\widehat{AB}:L^1_{h_1,\hat{h}}=\{(x_1,x_2)|H(x_1,x_2)=h_1,x_1\geq0,x_2\geq0,h_1>0\},\\
&\widehat{BC}:L^2_{h_1,\hat{h}}=\{(x_1,x_2)|H(x_1,x_2)=h_1,x_1>0,x_2<0,h_1>0\},\\
&\widehat{CE}:L^3_{h_1,\hat{h}}=\{(x_1,x_2)|H(x_1,x_2)=h_1,x_1<0,x_2<0,h_1>0\},\\
&\widehat{EA}:L^4_{h_1,\hat{h}}=\{(x_1,x_2)|H(x_1,x_2)=h_1,x_1<0,x_2>0,h_1>0\}.
\end{aligned}$$

Since $H(x_1,x_2)$ in (3.1) is independent of $\mathbf{y}$, the first order Melnikov function $M(h,\delta)$ of system (1.8) has the form of (2.19). We first calculate $M_1(h)=\sum\limits_{k=1}^4M_1^k(h)$, here $h=(h_1,\hat{h})^T=(h_1,h_3,h_4,\cdots,h_n)^T$. Put
$$I^i_{k_1k_2}(h_1)=\int_{L^i_{h_1,\hat{h}}}x_1^{k_1}x_2^{k_2}dx_2,\ i=1,2,3,4.$$
Let $O$ be the coordinate origin, by Green's formula two times, one has
\begin{eqnarray*}
\begin{aligned}
\int_{L^1_{h_1,\hat{h}}}x_1^{k_1}x_2^{k_2}dx_1=&\oint_{L^1_{h_1,\hat{h}}\cup \overrightarrow{BO}\cup \overrightarrow{OA}}x_1^{k_1}x_2^{k_2}dx_1-\int_{ \overrightarrow{BO}}x_1^{k_1}x_2^{k_2}dx_1\\
=&-k_2\iint_{\textup{int}(L^1_{h_1,\hat{h}}\cup \overrightarrow{BO}\cup \overrightarrow{OA})}x_1^{k_1}x_2^{k_2-1}dx_1dx_2-\int_{ \overrightarrow{BO}}x_1^{k_1}x_2^{k_2}dx_1,\\
\int_{L^1_{h_1,\hat{h}}}x_1^{k_1+1}x_2^{k_2-1}dx_2=&\oint_{L^1_{h_1,\hat{h}}\cup \overrightarrow{BO}\cup \overrightarrow{OA}}x_1^{k_1+1}x_2^{k_2-1}dx_2\\
=&(k_1+1)\iint_{\textup{int}(L^1_{h_1,\hat{h}}\cup \overrightarrow{BO}\cup \overrightarrow{OA})}x_1^{k_1}x_2^{k_2-1}dx_1dx_2,
\end{aligned}
\end{eqnarray*}
which follow directly
\begin{eqnarray}
\int_{L^1_{h_1,\hat{h}}}x_1^{k_1}x_2^{k_2}dx_1=\begin{cases}
-\frac{k_2}{k_1+1}I^1_{k_1+1,k_2-1}(h_1),\ k_1\geq0,k_2\geq1,\\[0.1cm]
\frac{1}{k_1+1}(2h_1)^\frac{k_1+1}{2},\qquad\quad k_1\geq0,k_2=0.\\
\end{cases}
\end{eqnarray}
Similarly, one can get
\begin{eqnarray}\begin{aligned}
&\int_{L^2_{h_1,\hat{h}}}x_1^{k_1}x_2^{k_2}dx_1=\begin{cases}
-\frac{k_2}{k_1+1}I^2_{k_1+1,k_2-1}(h_1),\ k_1\geq0,k_2\geq1,\\[0.1cm]
-\frac{1}{k_1+1}(2h_1)^\frac{k_1+1}{2},\qquad\ k_1\geq0,k_2=0,\\
\end{cases}\\
&\int_{L^3_{h_1,\hat{h}}}x_1^{k_1}x_2^{k_2}dx_1=\begin{cases}
-\frac{k_2}{k_1+1}I^3_{k_1+1,k_2-1}(h_1),\ k_1\geq0,k_2\geq1,\\[0.1cm]
-\frac{(-1)^{k_1}}{k_1+1}(2h_1)^\frac{k_1+1}{2},\quad\ \ k_1\geq0,k_2=0,\\
\end{cases}\\
&\int_{L^4_{h_1,\hat{h}}}x_1^{k_1}x_2^{k_2}dx_1=\begin{cases}
-\frac{k_2}{k_1+1}I^4_{k_1+1,k_2-1}(h_1),\ k_1\geq0,k_2\geq1,\\[0.1cm]
\frac{(-1)^{k_1}}{k_1+1}(2h_1)^\frac{k_1+1}{2},\qquad\ k_1\geq0,k_2=0.\\
\end{cases}\end{aligned}\end{eqnarray}
Thus, by (2.19), (3.4) and (3.5), one has
{\small\begin{eqnarray}
\begin{aligned}
&M_1(h)=\int_{L^1_{h_1,\hat{h}}}g_2^1(x_1,x_2,h_3,\cdots,h_n)dx_1-g_1^1(x_1,x_2,h_3,\cdots,h_n)dx_2\\
&+\int_{L^2_{h_1,\hat{h}}}g_2^2(x_1,x_2,h_3,\cdots,h_n)dx_1-g_1^2(x_1,x_2,h_3,\cdots,h_n)dx_2\\
&+\int_{L^3_{h_1,\hat{h}}}g_2^3(x_1,x_2,h_3,\cdots,h_n)dx_1-g_1^3(x_1,x_2,h_3,\cdots,h_n)dx_2\\
&+\int_{L^4_{h_1,\hat{h}}}g_2^4(x_1,x_2,h_3,\cdots,h_n)dx_1-g_1^4(x_1,x_2,h_3,\cdots,h_n)dx_2\\
=&\sum\limits_{\substack{k_1+k_2+\cdots+k_n=0\\ k_2\neq0}}^m\frac{k_2}{k_1+1} a^2_{k_1k_2\cdots k_n}h_3^{k_3}\cdots h_n^{k_n}I^1_{k_1+1,k_2-1}-
\sum\limits_{k_1+k_2+\cdots+k_n=0}^m a^1_{k_1k_2\cdots k_n}h_3^{k_3}\cdots h_n^{k_n}I^1_{k_1k_2}\\&+
\sum\limits_{k_1+k_3+\cdots+k_n=0}^m\frac{1}{k_1+1} a^2_{k_10k_3\cdots k_n}h_3^{k_3}\cdots h_n^{k_n}(2h_1)^\frac{k_1+1}{2}\\
&+\sum\limits_{\substack{k_1+k_2+\cdots+k_n=0\\ k_2\neq0}}^m\frac{k_2}{k_1+1} b^2_{k_1k_2\cdots k_n}h_3^{k_3}\cdots h_n^{k_n}I^2_{k_1+1,k_2-1}-
\sum\limits_{k_1+k_2+\cdots+k_n=0}^m b^1_{k_1k_2\cdots k_n}h_3^{k_3}\cdots h_n^{k_n}I^2_{k_1k_2}\\&-
\sum\limits_{k_1+k_3+\cdots+k_n=0}^m\frac{1}{k_1+1} b^2_{k_10k_3\cdots k_n}h_3^{k_3}\cdots h_n^{k_n}(2h_1)^\frac{k_1+1}{2}\\
&+\sum\limits_{\substack{k_1+k_2+\cdots+k_n=0\\ k_2\neq0}}^m\frac{k_2}{k_1+1} c^2_{k_1k_2\cdots k_n}h_3^{k_3}\cdots h_n^{k_n}I^3_{k_1+1,k_2-1}-
\sum\limits_{k_1+k_2+\cdots+k_n=0}^m c^1_{k_1k_2\cdots k_n}h_3^{k_3}\cdots h_n^{k_n}I^3_{k_1k_2}\\&-
\sum\limits_{k_1+k_3+\cdots+k_n=0}^m\frac{(-1)^{k_1}}{k_1+1} c^2_{k_10k_3\cdots k_n}h_3^{k_3}\cdots h_n^{k_n}(2h_1)^\frac{k_1+1}{2}\\
&+\sum\limits_{\substack{k_1+k_2+\cdots+k_n=0\\ k_2\neq0}}^m\frac{k_2}{k_1+1} d^2_{k_1k_2\cdots k_n}h_3^{k_3}\cdots h_n^{k_n}I^4_{k_1+1,k_2-1}-
\sum\limits_{k_1+k_2+\cdots+k_n=0}^m d^1_{k_1k_2\cdots k_n}h_3^{k_3}\cdots h_n^{k_n}I^4_{k_1k_2}\\&+
\sum\limits_{k_1+k_3+\cdots+k_n=0}^m\frac{(-1)^{k_1}}{k_1+1} d^2_{k_10k_3\cdots k_n}h_3^{k_3}\cdots h_n^{k_n}(2h_1)^\frac{k_1+1}{2}\\
:=&\sum\limits_{k_1+k_2+\cdots+k_n=0}^mh_3^{k_3}\cdots h_n^{k_n}\Big[\alpha^1_{k_1k_2\cdots k_n}I^1_{k_1,k_2}(h_1)+
\alpha^2_{k_1k_2\cdots k_n}I^2_{k_1,k_2}(h_1)+\alpha^3_{k_1k_2\cdots k_n}I^3_{k_1,k_2}(h_1)\\&+\alpha^4_{k_1k_2\cdots k_n}I^4_{k_1,k_2}(h_1)\Big]
+\sum\limits_{k_1+k_3+\cdots+k_n=0}^mh_3^{k_3}\cdots h_n^{k_n}\beta_{k_10k_3\cdots k_n}h_1^\frac{k_1+1}{2},
\end{aligned}
\end{eqnarray}}
where
$$\begin{aligned}
\alpha^1_{k_1k_2\cdots k_n}=&\begin{cases}
-a^1_{k_1k_2\cdots k_n}-\frac{k_2+1}{k_1}a^2_{k_1-1,k_2+1,k_3\cdots k_n},\ k_1\neq0,\\
-a^1_{k_1k_2\cdots k_n},\ \, \qquad\qquad\qquad\qquad\qquad k_1=0,\\
\end{cases}\\
\alpha^2_{k_1k_2\cdots k_n}=&\begin{cases}
-b^1_{k_1k_2\cdots k_n}-\frac{k_2+1}{k_1}b^2_{k_1-1,k_2+1,k_3\cdots k_n},\ k_1\neq0,\\
-b^1_{k_1k_2\cdots k_n},\ \, \qquad\qquad\qquad\qquad\qquad k_1=0,\\
\end{cases}\\
\alpha^3_{k_1k_2\cdots k_n}=&\begin{cases}
-c^1_{k_1k_2\cdots k_n}-\frac{k_2+1}{k_1}c^2_{k_1-1,k_2+1,k_3\cdots k_n},\ k_1\neq0,\\
-c^1_{k_1k_2\cdots k_n},\ \, \qquad\qquad\qquad\qquad\qquad k_1=0,\\
\end{cases}\\
\alpha^4_{k_1k_2\cdots k_n}=&\begin{cases}
-d^1_{k_1k_2\cdots k_n}-\frac{k_2+1}{k_1}d^2_{k_1-1,k_2+1,k_3\cdots k_n},\ k_1\neq0,\\
-d^1_{k_1k_2\cdots k_n},\ \, \qquad\qquad\qquad\qquad\qquad k_1=0,
\end{cases}\\
\beta_{k_10k_3\cdots k_n}=&\frac{1}{k_1+1}2^\frac{k_1+1}{2}\big[a^2_{k_10k_3\cdots k_n}-b^2_{k_10k_3\cdots k_n}-(-1)^{k_1}c^2_{k_10k_3\cdots k_n}\\&+(-1)^{k_1}d^2_{k_10k_3\cdots k_n}\big].
\end{aligned}$$
%$$\beta_{k_10k_3\cdots k_n}=\frac{1}{k_1+1}2^\frac{k_1+1}{2}\big(a^2_{k_10k_3\cdots k_n}-b^2_{k_10k_3\cdots k_n}-(-1)^{k_1}c^2_{k_10k_3\cdots k_n}+(-1)^{k_1}d^2_{k_10k_3\cdots k_n}\big).$$
It is easy to see that $\alpha^i_{k_1k_2\cdots k_n}(i=1,2,3,4)$ and $\beta_{k_10k_3\cdots k_n}$  can be chosen arbitrarily.
\vskip 0.2 true cm

\noindent
{\bf Lemma 3.1.}\, {\it \begin{eqnarray}\begin{aligned}
&I^1_{k_1k_2}(h_1)=\begin{cases}
\gamma^{11}_{k_1k_2}h_1^\frac{k_1+k_2}{2}I^1_{00}(h_1)+\gamma^{12}_{k_1k_2}h_1^\frac{k_1+k_2-2}{2}I^1_{11}(h_1),\quad k_1+k_2\ \textup{even},\\[0.1cm]
\gamma^{13}_{k_1k_2}h_1^\frac{k_1+k_2-1}{2}I^1_{01}(h_1)+\gamma^{14}_{k_1k_2}h_1^\frac{k_1+k_2-1}{2}I^1_{10}(h_1),\ k_1+k_2\ \textup{odd},\\
\end{cases}\\[2mm]
&I^2_{k_1k_2}(h_1)=\begin{cases}
\gamma^{21}_{k_1k_2}h_1^\frac{k_1+k_2}{2}I^2_{00}(h_1)+\gamma^{22}_{k_1k_2}h_1^\frac{k_1+k_2-2}{2}I^2_{11}(h_1),\quad k_1+k_2\ \textup{even},\\[0.1cm]
\gamma^{23}_{k_1k_2}h_1^\frac{k_1+k_2-1}{2}I^2_{01}(h_1)+\gamma^{24}_{k_1k_2}h_1^\frac{k_1+k_2-1}{2}I^2_{10}(h_1),\ k_1+k_2\ \textup{odd},\\
\end{cases}\\[2mm]
&I^3_{k_1k_2}(h_1)=\begin{cases}
\gamma^{31}_{k_1k_2}h_1^\frac{k_1+k_2}{2}I^3_{00}(h_1)+\gamma^{32}_{k_1k_2}h_1^\frac{k_1+k_2-2}{2}I^3_{11}(h_1),\quad k_1+k_2\ \textup{even},\\[0.1cm]
\gamma^{33}_{k_1k_2}h_1^\frac{k_1+k_2-1}{2}I^3_{01}(h_1)+\gamma^{34}_{k_1k_2}h_1^\frac{k_1+k_2-1}{2}I^3_{10}(h_1),\ k_1+k_2\ \textup{odd},\\
\end{cases}\\[2mm]
&I^4_{k_1k_2}(h_1)=\begin{cases}
\gamma^{41}_{k_1k_2}h_1^\frac{k_1+k_2}{2}I^4_{00}(h_1)+\gamma^{42}_{k_1k_2}h_1^\frac{k_1+k_2-2}{2}I^4_{11}(h_1),\quad k_1+k_2\ \textup{even},\\[0.1cm]
\gamma^{43}_{k_1k_2}h_1^\frac{k_1+k_2-1}{2}I^4_{01}(h_1)+\gamma^{44}_{k_1k_2}h_1^\frac{k_1+k_2-1}{2}I^4_{10}(h_1),\ k_1+k_2\ \textup{odd},\\
\end{cases}\\
\end{aligned}
\end{eqnarray}
where $\gamma^{ij}_{k_1k_2}(i,j=1,2,3,4)$ are constants.}
\vskip 0.2 true cm

\noindent
{\bf Proof.} We only prove the first equality in (3.7). The others can be shown similarly. Differentiating both sides of the following equation
\begin{eqnarray}
H(x_1,x_2)=\frac{1}{2}(x_1^2+x_2^2)=h_1
\end{eqnarray}
with respect to $x_2$, one obtains
\begin{eqnarray}
x_2+x_1\frac{\partial x_1}{\partial x_2}=0.
\end{eqnarray}
Multiplying (3.9) by $x_1^{k_1}x_2^{k_2-1}dx_2$, integrating over $L^1_{h_1,\hat{h}}$ and noting that (3.4), one has
\begin{eqnarray}
I^1_{k_1k_2}(h_1)=\frac{k_2-1}{k_1+2}I^1_{k_1+2,k_2-2}(h_1), \ k_1\geq0,k_2\geq2.
\end{eqnarray}
Similarly, multiplying (3.8) by $x_1^{k_1-2}x_2^{k_2}dx_2$ and integrating over $L^1_{h_1,\hat{h}}$ give
\begin{eqnarray}
I^1_{k_1k_2}(h_1)=2h_1I^1_{k_1-2,k_2}(h_1)-I^1_{k_1-2,k_2+2}(h_1),\ k_1\geq2,k_2\geq0.
\end{eqnarray}
Eliminating $I^1_{k_1+2,k_2-2}(h_1)$ and $I^1_{k_1-2,k_2+2}(h_1)$ in (3.10) and (3.11) yields
\begin{eqnarray}
I^1_{k_1k_2}(h_1)=\frac{2k_1}{k_1+k_2+1}h_1I^1_{k_1-2,k_2}(h_1),\ k_1\geq2,k_2\geq0
\end{eqnarray}
and
\begin{eqnarray}
I^1_{k_1k_2}(h_1)=\frac{2(k_2-1)}{k_1+k_2+1}h_1I^1_{k_1,k_2-2}(h_1),\ k_1\geq0,k_2\geq2.
\end{eqnarray}

Now we use mathematical induction to prove the first equality in (3.7). From (3.12) and (3.13), one can get
\begin{eqnarray}
\begin{cases}
I^1_{0,2}(h_1)=\frac{2}{3}h_1I^1_{0,0}(h_1),\ \ I^1_{2,0}(h_1)=\frac{4}{3}h_1I^1_{0,0}(h_1),\\
I^1_{0,3}(h_1)=h_1I^1_{0,1}(h_1),\ \quad I^1_{1,2}(h_1)=\frac{1}{2}h_1I^1_{1,0}(h_1),\\
I^1_{2,1}(h_1)=h_1I^1_{0,1}(h_1),\ \quad I^1_{3,0}(h_1)=\frac{3}{2}h_1I^1_{1,0}(h_1),
\end{cases}
\end{eqnarray}
which follow directly that the first equality in (3.7) holds for $k_1+k_2=2,3$. Suppose that the first equality in (3.7) holds when $k_1+k_2\leq l-1$, where $l\geq4$ is an even number. Then, taking $(k_1,k_2)=(0,l),(1,l-1),(2,l-2),\cdots,(l-2,2)$ in (3.13) and $(k_1,k_2)=(l-1,1),(l,0)$ in (3.12), one gets
\begin{eqnarray}
\left(\begin{matrix}
                 I^1_{0,l}(h_1)\\
                 I^1_{1,l-1}(h_1)\\
                 I^1_{2,l-2}(h_1)\\
                 \vdots\\
                  I^1_{l-1,1}(h_1)\\
                  I^1_{l,0}(h_1)
                   \end{matrix}\right)\ \
=\frac{2}{l+1}\left(\begin{matrix}
                (l-1)h_1I^1_{0,l-2}(h_1)\\
                (l-2)h_1I^1_{1,l-3}(h_1)\\
                 (l-3)h_1I^1_{2,l-4}(h_1)\\
                  \vdots\\
                   (l-1)h_1I^1_{l-3,1}(h_1)\\
                    lh_1I^1_{l-2,0}(h_1)
                \end{matrix}\right).
\end{eqnarray}
Hence, by (3.15), one has for $k_1+k_2=l$
\begin{eqnarray*}
\begin{aligned}
I_{k_1,k_2}(h_1)=&h_1[\gamma^1_{k_1k_2}h_1^{\frac{l-2}{2}}I_{0,0}(h_1)+\gamma^2_{k_1k_2}h_1^{\frac{l-4}{2}}I_{0,0}(h_1)]\\
=&\gamma^{11}_{k_1k_2}h_1^{\frac{l}{2}}I_{0,0}(h_1)+\gamma^{12}_{k_1k_2}h_1^{\frac{l-2}{2}}I_{0,0}(h_1).
\end{aligned}
\end{eqnarray*}

Similar to the above proof, if $k_1+k_2=l$ is an odd number, we can prove that the second equality in (3.7) holds too. This ends the proof. \quad $\lozenge$
\vskip 0.2 true cm

\noindent
{\bf Lemma 3.2.}\, {\it $M_1(h)$ can be expressed as
\begin{eqnarray}
M_1(h)=\sum\limits_{k_1+k_2+\cdots+k_n=0}^m\lambda^1_{k_1k_2\cdots k_n}h^{\frac{k_1+k_2+1}{2}}h_3^{k_3}\cdots h_n^{k_n},
\end{eqnarray}
here $\lambda^1_{k_1k_2\cdots k_n}$ is a constant.}
\vskip 0.2 true cm

\noindent
{\bf Proof.} By (3.6) and Lemma 3.1, one has
\begin{eqnarray}\begin{aligned}
M_1(h)=&\sum\limits_{\substack{k_1+k_2+\cdots+k_n=0\\ k_1+k_2=0 mod 2}}^m h_3^{k_3}\cdots h_n^{k_n}\Big[\alpha^1_{k_1k_2\cdots k_n}\Big(\gamma^{11}_{k_1k_2}h_1^\frac{k_1+k_2}{2}I^1_{00}+\gamma^{12}_{k_1k_2}h_1^\frac{k_1+k_2-2}{2}I^1_{11}\Big)\\
&+\alpha^2_{k_1k_2\cdots k_n}\Big(\gamma^{21}_{k_1k_2}h_1^\frac{k_1+k_2}{2}I^2_{00}+\gamma^{22}_{k_1k_2}h_1^\frac{k_1+k_2-2}{2}I^2_{11}\Big)\\
&+\alpha^3_{k_1k_2\cdots k_n}\Big(\gamma^{31}_{k_1k_2}h_1^\frac{k_1+k_2}{2}I^3_{00}+\gamma^{32}_{k_1k_2}h_1^\frac{k_1+k_2-2}{2}I^3_{11}\Big)\\
&+\alpha^4_{k_1k_2\cdots k_n}\Big(\gamma^{41}_{k_1k_2}h_1^\frac{k_1+k_2}{2}I^4_{00}+\gamma^{42}_{k_1k_2}h_1^\frac{k_1+k_2-2}{2}I^4_{11}\Big)\Big]\\
&+\sum\limits_{\substack{k_1+k_2+\cdots+k_n=0\\ k_1+k_2=1 mod 2}}^m h_3^{k_3}\cdots h_n^{k_n}\Big[\alpha^1_{k_1k_2\cdots k_n}\Big(\gamma^{13}_{k_1k_2}h_1^\frac{k_1+k_2-1}{2}I^1_{01}+\gamma^{14}_{k_1k_2}h_1^\frac{k_1+k_2-1}{2}I^1_{10}\Big)\\
&+\alpha^2_{k_1k_2\cdots k_n}\Big(\gamma^{23}_{k_1k_2}h_1^\frac{k_1+k_2-1}{2}I^2_{01}+\gamma^{24}_{k_1k_2}h_1^\frac{k_1+k_2-1}{2}I^2_{10}\Big)\\
&+\alpha^3_{k_1k_2\cdots k_n}\Big(\gamma^{33}_{k_1k_2}h_1^\frac{k_1+k_2-1}{2}I^3_{01}+\gamma^{34}_{k_1k_2}h_1^\frac{k_1+k_2-1}{2}I^3_{10}\Big)\\
&+\alpha^4_{k_1k_2\cdots k_n}\Big(\gamma^{43}_{k_1k_2}h_1^\frac{k_1+k_2-1}{2}I^4_{01}+\gamma^{44}_{k_1k_2}h_1^\frac{k_1+k_2-1}{2}I^4_{10}\Big)\Big].
\end{aligned}\end{eqnarray}
By a straightforward calculation, we have
$$\begin{aligned}
&I^1_{00}=-\sqrt{2h_1},\ I^1_{10}=-\frac{\pi}{2}h_1,\  I^1_{01}=-h_1,\  I^1_{11}=-\frac{2}{3}\sqrt{2}h_1^\frac{3}{2},\\
&I^2_{00}=-\sqrt{2h_1},\ I^2_{10}=-\frac{\pi}{2}h_1,\  I^2_{01}=h_1,\  I^2_{11}=\frac{2}{3}\sqrt{2}h_1^\frac{3}{2},\\
&I^3_{00}=\sqrt{2h_1},\ I^3_{10}=-\frac{\pi}{2}h_1,\  I^3_{01}=-h_1,\  I^3_{11}=\frac{2}{3}\sqrt{2}h_1^\frac{3}{2},\\
&I^4_{00}=\sqrt{2h_1},\ I^4_{10}=-\frac{\pi}{2}h_1,\  I^4_{01}=h_1,\  I^4_{11}=-\frac{2}{3}\sqrt{2}h_1^\frac{3}{2}.
\end{aligned}$$
Therefore, (3.16) follows from (3.17) and the above equalities. This ends the proof.\quad $\lozenge$
\vskip 0.2 true cm

Now, let's go back to calculate $M_i(h)=\sum\limits_{k=1}^4M_i^k(h)(i=2,3,\cdots,n-1)$ in (2.19). Their detailed expressions are given in the following lemma.
\vskip 0.2 true cm

\noindent
{\bf Lemma 3.3.}\, {\it $M_i(h)(i=2,3,\cdots,n-1)$ can be written as
\begin{eqnarray}
\begin{aligned}
&M_2(h)=\sum\limits_{k_1+k_2+\cdots+k_n=0}^m\lambda^2_{k_1k_2\cdots k_n}h^{\frac{k_1+k_2}{2}}h_3^{k_3}\cdots h_n^{k_n},\\
&M_3(h)=\sum\limits_{k_1+k_2+\cdots+k_n=0}^m\lambda^3_{k_1k_2\cdots k_n}h^{\frac{k_1+k_2}{2}}h_3^{k_3}\cdots h_n^{k_n},\\
&\cdots\\
&M_{n-1}(h)=\sum\limits_{k_1+k_2+\cdots+k_n=0}^m\lambda^{n-1}_{k_1k_2\cdots k_n}h^{\frac{k_1+k_2}{2}}h_3^{k_3}\cdots h_n^{k_n},\\
\end{aligned}
\end{eqnarray}
where $\lambda^{i}_{k_1k_2\cdots k_n}(i=2,3,\cdots,n-1)$ are constants.}
\vskip 0.2 true cm

\noindent
{\bf Proof.} Without loss of generality, we only calculate $M_2(h)$ and the others can be calculated similarly. For the sake of clarity, we use the notations given in Lemma 2.1. So we first compute $\bar{R}^1_{3}(x_1,x_2,\hat{h})$ and $\bar{M}^1_2(h)$. It is easy to get that
$$\begin{aligned}
\bar{R}^1_{3}(x_1,x_2,\hat{h})=&\int_0^{x_2}g^1_{3}(x_1,x_2,\hat{h})dx_2\\
&=\sum\limits_{k_1+k_2+\cdots+k_n=0}^ma^3_{k_1k_2\cdots k_n}h^{k_3}_3\cdots h^{k_n}_n\int_0^{x_2}x_1^{k_1}x_2^{k_2}dx_2\\
&=\sum\limits_{k_1+k_2+\cdots+k_n=0}^m\frac{1}{k_2+1}a^3_{k_1k_2\cdots k_n}h^{k_3}_3\cdots h^{k_n}_nx_1^{k_1}x_2^{k_2+1}.
\end{aligned}$$
By (3.4) and Lemma 3.1, one has
$$\begin{aligned}
\bar{M}^1_2(h)=&\int_{L^1_{h_1,\hat{h}}}\bar{R}^1_{3}(x_1,x_2,\hat{h})dx_1\\
=&\sum\limits_{k_1+k_2+\cdots+k_n=0}^m\frac{1}{k_2+1}a^3_{k_1k_2\cdots k_n}h^{k_3}_3\cdots h^{k_n}_n\int_{L^1_{h_1,\hat{h}}}x_1^{k_1}x_2^{k_2+1}dx_1\\
=&-\sum\limits_{k_1+k_2+\cdots+k_n=0}^m\frac{1}{k_1+1}a^3_{k_1k_2\cdots k_n}h^{k_3}_3\cdots h_n^{k_n}I^1_{k_1+1,k_2}(h_1)\\
=&\sum\limits_{k_1+k_2+\cdots+k_n=0}^m\bar{\lambda}^2_{k_1k_2\cdots k_n}h_1^{\frac{k_1+k_2+2}{2}}h^{k_3}_3\cdots h_n^{k_n},
\end{aligned}$$
where $\bar{\lambda}^2_{k_1k_2\cdots k_n}$ is a constant.
Thus, by Lemma 2.1, it follows that
\begin{eqnarray}M^1_2(h)=\frac{\partial \bar{M}_2^1(h)}{\partial h_1}=\sum\limits_{k_1+k_2+\cdots+k_n=0}^m\bar{\lambda}^2_{k_1k_2\cdots k_n}h_1^{\frac{k_1+k_2}{2}}h^{k_3}_3\cdots h_n^{k_n}.\end{eqnarray}
Similarly, one obtains
\begin{eqnarray}\begin{aligned}
&M^2_2(h)=\sum\limits_{k_1+k_2+\cdots+k_n=0}^m\tilde{\lambda}^2_{k_1k_2\cdots k_n}h_1^{\frac{k_1+k_2}{2}}h^{k_3}_3\cdots h_n^{k_n},\\
&M^3_2(h)=\sum\limits_{k_1+k_2+\cdots+k_n=0}^m\hat{\lambda}^2_{k_1k_2\cdots k_n}h_1^{\frac{k_1+k_2}{2}}h^{k_3}_3\cdots h_n^{k_n},\\
&M^4_2(h)=\sum\limits_{k_1+k_2+\cdots+k_n=0}^m\check{\lambda}^2_{k_1k_2\cdots k_n}h_1^{\frac{k_1+k_2}{2}}h^{k_3}_3\cdots h_n^{k_n},\end{aligned}\end{eqnarray}
where $\tilde{\lambda}^2_{k_1k_2\cdots k_n}$, $\hat{\lambda}^2_{k_1k_2\cdots k_n}$ and $\check{\lambda}^2_{k_1k_2\cdots k_n}$ are constants.
Therefore, $M_2(h)$ in (3.18) follows from (2.19), (3.19) and (3.20). This completes the proof.\quad $\lozenge$
\vskip 0.2 true cm

In order to prove Theorem 1.2, we need the following result which was given by Han, Sun and Balanov in \cite{HSB}.
\vskip 0.2 true cm

\noindent
{\bf Lemma 3.4.} {\it Let $f_1,f_2,\cdots,f_n\in \mathbb{R}[x_1,x_2,\cdots,x_n]$ be real polynomials and assume that the map $g=(f_1,f_2,\cdots,f_n):\mathbb{R}^n\rightarrow\mathbb{R}^n$ admits finitely many zeros. Then, the number of zeros of $f$ is at most $\deg f_1\times \deg f_2\times\cdots\times\deg f_n$.}
\vskip 0.2 true cm

\noindent
{\bf Proof of Theorem 1.2.} If $m\geq2$, by Lemmas 3.2 and 3.3, one gets that the first order Melnikov vector function $M(h)$ of system (1.8) has the form
\begin{eqnarray}
M(h)=\left(\begin{matrix}
                M_1(h)\\
                 M_2(h)\\
                 M_3(h)\\
                 \vdots\\
                  M_{n-1}(h)
                   \end{matrix}\right)=\left(\begin{matrix}
                 h_2\sum\limits_{k_1+k_2+\cdots+k_n=0}^m\lambda^1_{k_1k_2\cdots k_n}h_2^{k_1+k_2}h_3^{k_3}\cdots h_n^{k_n}\\
                \sum\limits_{k_1+k_2+\cdots+k_n=0}^m{\lambda}^2_{k_1k_2\cdots k_n}h_2^{k_1+k_2}h_3^{k_3}\cdots h_n^{k_n}\\
                 \sum\limits_{k_1+k_2+\cdots+k_n=0}^m{\lambda}^3_{k_1k_2\cdots k_n}h_2^{k_1+k_2}h_3^{k_3}\cdots h_n^{k_n}\\
                 \vdots\\
                 \sum\limits_{k_1+k_2+\cdots+k_n=0}^m{\lambda}^{n-1}_{k_1k_2\cdots k_n}h_2^{k_1+k_2}h_3^{k_3}\cdots h_n^{k_n}
\end{matrix}\right),
\end{eqnarray}
where $h_2=\sqrt{h_1}>0$. Let $M_1(h)=h_2\tilde{M}_1(h)$. Then $\tilde{M}_1(h)$ is a polynomial of $h=(h_2,h_3,\cdots,h_n)^T$ with degree $m$ and so are $M_k(h)(k=2,3,\cdots,n-1)$. By Lemma 3.4, one has that $M(h)$ has at most $m^{n-1}$ zeros. Hence, the conclusion follows from Theorem 1.1.

If $m=1$, for the convenience of writing, we introduce a shorthand notation and rewrite system (1.8) as follows
\begin{eqnarray}\begin{aligned}
&\begin{cases}
\dot{x}_1=x_2+\varepsilon (a^1_0+a^1_1x_1+\cdots+a^1_nx^n),\\
\dot{x}_2=-x_1+\varepsilon (a^2_0+a^2_1x_1+\cdots+a^2_nx^n),\\
\dot{x}_3=\varepsilon (a^3_0+a^3_1x_1+\cdots+a^3_nx^n),\\
\cdots\\
\dot{x}_n=\varepsilon (a^n_0+a^n_1x_1+\cdots+a^n_nx^n),\\
\end{cases}x_1\geq0,x_2\geq0,\\
&\begin{cases}
\dot{x}_1=x_2+\varepsilon (b^1_0+b^1_1x_1+\cdots+b^1_nx^n),\\
\dot{x}_2=-x_1+\varepsilon (b^2_0+b^2_1x_1+\cdots+b^2_nx^n),\\
\dot{x}_3=\varepsilon (b^3_0+b^3_1x_1+\cdots+b^3_nx^n),\\
\cdots\\
\dot{x}_n=\varepsilon (b^n_0+b^n_1x_1+\cdots+b^n_nx^n),\\
\end{cases}x_1>0,x_2<0,\\
&\begin{cases}
\dot{x}_1=x_2+\varepsilon (c^1_0+c^1_1x_1+\cdots+c^1_nx^n),\\
\dot{x}_2=-x_1+\varepsilon (c^2_0+c^2_1x_1+\cdots+c^2_nx^n),\\
\dot{x}_3=\varepsilon (c^3_0+c^3_1x_1+\cdots+c^3_nx^n),\\
\cdots\\
\dot{x}_n=\varepsilon (c^n_0+c^n_1x_1+\cdots+c^n_nx^n),\\
\end{cases}x_1<0,x_2<0,\\
&\begin{cases}
\dot{x}_1=x_2+\varepsilon (d^1_0+d^1_1x_1+\cdots+d^1_nx^n),\\
\dot{x}_2=-x_1+\varepsilon (d^2_0+d^2_1x_1+\cdots+d^2_nx^n),\\
\dot{x}_3=\varepsilon (d^3_0+d^3_1x_1+\cdots+d^3_nx^n),\\
\cdots\\
\dot{x}_n=\varepsilon (d^n_0+d^n_1x_1+\cdots+d^n_nx^n),\\
\end{cases}x_1<0,x_2>0.
\end{aligned}\end{eqnarray}
Similar to the above calculation, one can obtain the first order Melnikov function of system (3.22) as follows
\begin{eqnarray}
\begin{aligned}
M(h)=\big(M_1(h),M_2(h),\cdots, M_{n-1}(h) \big)^T,
\end{aligned}
\end{eqnarray}
where
\begin{eqnarray*}
\begin{aligned}
&M_1(h)=\sqrt{h_1}\big(\rho_1^1+\rho_2^1\sqrt{h_1}+\rho_3^1 h_3+\cdots+\rho_n^1 h_n\big),\\
&M_2(h)=\rho_1^2+\rho_2^2\sqrt{h_1}+\rho_3^2 h_3+\cdots+\rho_n^2 h_n,\\
&\cdots\\
&M_{n-1}(h)=\rho_1^{n-1}+\rho_2^{n-1}\sqrt{h_1}+\rho_3^{n-1} h_3+\cdots+\rho_n^{n-1} h_n,\\
\end{aligned}
\end{eqnarray*}
\begin{eqnarray*}
\begin{aligned}
\rho_1^1=&a_0^2+a_0^1-b_0^2+b_0^1-c_0^2-c_0^2+d_0^2-d_0^1,\\
\rho_2^1=&a_1^2-b_1^2+c_1^2-d_1^2+a_2^1-b_2^1+c_2^1-d_2^1\\&+\frac{\pi}{2}(a_2^2-b_2^2+c_2^2-d_2^2+a_1^1+b_1^1+c_1^1+d_1^1),\\
\rho_k^1=&\sqrt{2}(a_k^2+a_k^1-b_k^2+b_k^1-c_k^2-c_k^1+d_k^2-d_k^1),\\
\rho_1^i=&\frac{\pi}{2}(a_0^{i+1}-b_0^{i+1}+c_0^i-d_0^{i+1}),\\
\rho_2^i=&\sqrt{2}(a_1^{i+1}+a_2^{i+1}+b_1^{i+1}-b_2^{i+1}-c_1^{i+1}-c_2^{i+1}-d_1^{i+1}+d_2^{i+1}),\\
\rho_k^i=&\frac{\pi}{2}(a_k^{i+1}-b_k^{i+1}+c_k^{i+1}-d_k^{i+1}),\\
k=&3,4,\cdots,n;i=2,3,\cdots,n-1.
\end{aligned}
\end{eqnarray*}
It is easy to check that the coefficients $\rho_j^i$ of $M_i(h)$ in (3.23) can be chosen arbitrarily for $i=1,2,\cdots,n-1;j=1,2,\cdots,n$. Let $a^i_j$, $b^i_j$, $c^i_j$ and $d^i_j$ in (3.22) be zero except $a_0^1>0$, $a_2^2<0$, $a_0^k$ and $a_k^k(k=3,\cdots,n)$. By Cramer's Rule, one can get that linear system of equations $M(h)=0$ has a unique root $$X_0=(\sqrt{h_1},h_3,\cdots,h_n)=\Big(-\frac{a_0^1}{a^2_2},-\frac{a_0^3}{a^3_3},\cdots,-\frac{a_0^n}{a^n_n}\Big).$$ Further, one has
$$\det\frac{\partial \big(M_1(h),M_2(h),\cdots,M_{n-1}(h)\big)}{\partial(h_1,h_3,\cdots,h_n)}\Big| _{X_0}=-\Big(\frac{\pi}{2}\Big)^{n-1}\frac{a_2^2}{2a_0^1}\prod\limits_{k=2}^na^k_k\neq0.$$
Then, Theorem 1.1 gives the desired result. This ends the proof of Theorem 1.2. \quad $\lozenge$

\vskip 0.2 true cm

\noindent
{\bf Acknowledgment}
 \vskip 0.2 true cm

\noindent
The authors would like to express their sincere appreciation to the referee for his/her valuable suggestions and comments. This work was supported by National Natural Science Foundation of China(11701306), Construction of First-class Disciplines of Higher Education of Ningxia(Pedagogy)(NXYLXK2017B11), Ningxia Natural Science Foundation(2019AAC03247), Young Top-notch Talent of Ningxia and Key Program of Ningxia Normal University(NXSFZDA1901).

\newpage
\appendix
\noindent
\section{ Appendix}
\vskip 0.2 true cm
 \setcounter{equation}{0}
\renewcommand\theequation{A.\arabic{equation}}

Consider a piecewise smooth near-Hamiltonian system of the form
\begin{eqnarray}
\begin{cases}
\dot{x}=H_y(x,y)+\varepsilon p(x,y,\delta),\\
\dot{y}=-H_x(x,y)+\varepsilon q(x,y,\delta),\\
\end{cases}
\end{eqnarray}
where $0<|\varepsilon|\ll1$,
\begin{eqnarray*}
\big(H(x,y),p(x,y,\delta),q(x,y,\delta)\big)=\begin{cases}
\big(H^1(x,y),p^1(x,y,\delta),q^1(x,y,\delta)\big),\ x\geq0,y\geq0,\\
\big(H^2(x,y),p^2(x,y,\delta),q^2(x,y,\delta)\big),\ x>0,y<0,\\
\big(H^3(x,y),p^3(x,y,\delta),q^3(x,y,\delta)\big),\ x<0,y<0,\\
\big(H^4(x,y),p^4(x,y,\delta),q^4(x,y,\delta)\big),\ x<0,y>0\\
\end{cases}
\end{eqnarray*}
and  $H^i(x,y),p^i(x,y,\delta),q^i(x,y,\delta)\in C^\infty$, $i=1,2,3,4$, $\delta\in D\subset \mathbb{R}^m$ with $D$ compact. For an open interval $\Sigma$, we suppose that system (A.1) has a family of closed orbits $\{L_h |h\in\Sigma\}$ around the origin with a clockwise orientation and each closed orbit $L_h$ intersects $x$-axis (resp. $y$-axis) at two different points $A_1(h) = (0,a_1(h))$ and $A_3(h) = (0,a_3(h))$ (resp. $A_2(h) = (a_2(h),0)$ and $A_4(h) = (a_4(h),0))$ with $a_1(h) > 0 > a_3(h)$ and $a_2(h) > 0 > a_3(h)$, see Fig. 3.

Then, from \cite{WHC,HD}, one knows that the first order Melnikov function of system (A.1) is
\begin{eqnarray*}
\begin{aligned}
M(h)=&\int_{\widehat{A_1A_2}}q^1dx-p^1dy+\frac{H^1_y(A_1)H^3_y(A_3)H^4_x(A_4)}{H^4_y(A_1)H^2_y(A_3)H^3_x(A_4)}
\int_{\widehat{A_2A_3}}q^2dx-p^2dy\\
&+\frac{H^1_y(A_1)H^4_x(A_4)}{H^4_y(A_1)H^3_x(A_4)}
\int_{\widehat{A_3A_4}}q^3dx-p^3dy+\frac{H^1_y(A_1)}{H^4_y(A_1)}
\int_{\widehat{A_4A_1}}q^4dx-p^4dy,\ h\in\Sigma,
\end{aligned}
\end{eqnarray*}
where
$$\begin{aligned}
&\widehat{A_1A_2}=\{(x,y)|H^1(x,y)=h,x\geq0,y\geq0\},\\
&\widehat{A_2A_3}=\{(x,y)|H^2(x,y)=H^2(A_2),x>0,y<0\},\\
&\widehat{A_3A_4}=\{(x,y)|H^3(x,y)=H^3(A_3),x<0,y<0\},\\
&\widehat{A_4A_1}=\{(x,y)|H^4(x,y)=H^4(A_4),x<0,y>0\}.
\end{aligned}$$
Put
\begin{eqnarray}\begin{aligned}
&M_1(h,\delta)=\int_{\widehat{A_1A_2}}q^1(x,y,\delta)dx-p^1(x,y,\delta)dy,\\
&M_2(h,\delta)=\int_{\widehat{A_2A_3}}q^2(x,y,\delta)dx-p^2(x,y,\delta)dy,\\
&M_3(h,\delta)=\int_{\widehat{A_3A_4}}q^3(x,y,\delta)dx-p^3(x,y,\delta)dy,\\
&M_4(h,\delta)=\int_{\widehat{A_4A_1}}q^4(x,y,\delta)dx-p^4(x,y,\delta)dy.
\end{aligned}\end{eqnarray}
Then, one has the following lemma for proving the Lemma 2.1 in Section 2.
\vskip 0.2true cm

\noindent
{\bf Lemma A.1.} For the integrals $M_k(h,\delta)(k=1,2,3,4)$ in (A.2), we have
$$\begin{aligned}
&\frac{\partial M_1(h,\delta)}{\partial h}=\int_{\widehat{A_1A_2}}(p^1_x+q^1_y)dt+p^1(0,a_1(h),\delta)a'_1(h)+q^1(a_2(h),0,\delta)a'_2(h),\\
&\frac{\partial M_2(h,\delta)}{\partial h}=\int_{\widehat{A_2A_3}}(p^2_x+q^2_y)dtH^2_y(A_2)a'_2(h)-p^2(0,a_3(h),\delta)a'_3(h)-q^2(a_2(h),0,\delta)a'_2(h),\\
&\frac{\partial M_3(h,\delta)}{\partial h}=\int_{\widehat{A_3A_4}}(p^3_x+q^3_y)dtH^3_y(A_3)a'_3(h)+p^3(0,a_3(h),\delta)a'_3(h)+q^3(a_4(h),0,\delta)a'_4(h),\\
&\frac{\partial M_4(h,\delta)}{\partial h}=\int_{\widehat{A_4A_1}}(p^4_x+q^4_y)dtH^4_y(A_4)a'_4(h)-p^4(0,a_1(h),\delta)a'_1(h)-q^4(a_4(h),0,\delta)a'_4(h).
\end{aligned}$$

\noindent
{\bf Proof.} We only prove the first equality and the others can be proved similarly. By Green's formula, one has
$$\begin{aligned}
M_1(h,\delta)=&\oint_{\widehat{A_1A_2}\cup \overrightarrow{A_2O}\cup\overrightarrow{OA_1}}q^1(x,y,\delta)dx-p^1(x,y,\delta)dy\\&-\int_{\overrightarrow{A_2O}}q^1(x,0,\delta)dx
+\int_{\overrightarrow{OA_1}}p^1(0,y,\delta)dy,\\
=&\iint_{\textup{int}(\widehat{A_1A_2}\cup \overrightarrow{A_2O}\cup\overrightarrow{OA_1})}\big[p^1_x(x,y,\delta)+q^1_y(x,y,\delta)\big]dxdy
\\&-\int_{\overrightarrow{A_2O}}q^1(x,0,\delta)dx+\int_{\overrightarrow{OA_1}}p^1(0,y,\delta)dy,\\
=&\oint_{\widehat{A_1A_2}\cup \overrightarrow{A_2O}\cup\overrightarrow{OA_1}}\tilde{q}^1(x,y,\delta)dx-\int_{\overrightarrow{A_2O}}q^1(x,0,\delta)dx
+\int_{\overrightarrow{OA_1}}p^1(0,y,\delta)dy,\\
\end{aligned}$$
where $$\tilde{q}^1(x,y,\delta)=\int_0^y\big[p^1_x(x,s,\delta)+q^1_y(x,s,\delta)\big]ds.$$
In view of $$\int_{\small\overrightarrow{A_2O}}\tilde{q}^1(x,y,\delta)dx=\int_{\small\overrightarrow{OA_1}}\tilde{q}^1(x,y,\delta)dy=0,$$
one obtains
$$M_1(h,\delta)=\oint_{\widehat{A_1A_2}}\tilde{q}^1(x,y,\delta)dx-\int_{\overrightarrow{A_2O}}q^1(x,0,\delta)dx
+\int_{\overrightarrow{OA_1}}p^1(0,y,\delta)dy.$$

Without loss of generality, suppose that there exists only one point $B(h)=(b(h),\tilde{b}(h))$ on the curve $\widehat{A_1A_2}$ such that
$H^1_y(B(h))\neq0,\ H^1_x(B(h))=0$. Then, the curve $\widehat{A_1B}$ can be expressed as $y=y_1(x,h)$ and the curve $\widehat{BA_2}$ can be expressed as $y=y_2(x,h)$ for $0<x<c(h)$, see Fig. 4.
Then, one has
$$\int_{\widehat{A_1A_2}}\tilde{q}^1(x,y,\delta)dx=\int_0^{c(h)}\tilde{q}^1(x,y_1(x,h),\delta)dx+
\int^{a_2(h)}_{c(h)}\tilde{q}^1(x,y_2(x,h),\delta)dx,$$
which implies
$$\begin{aligned}\frac{d}{dh}\Big(\int_{\widehat{A_1A_2}}\tilde{q}^1(x,y,\delta)dx\Big)=&\int_0^{c(h)}\tilde{q}^1_y(x,y_1(x,h),\delta)\frac{\partial y_1}{\partial h}dx+
\int^{a_2(h)}_{c(h)}\tilde{q}^1_y(x,y_2(x,h),\delta)\frac{\partial y_2}{\partial h}dx\\
=&\int_{\widehat{A_1A_2}}\tilde{q}^1_y(x,y,\delta)\frac{\partial y}{\partial h}dx.\end{aligned}$$
 On the other hand, along $\widehat{A_1A_2}$, one gets $H^1_y\frac{\partial y}{\partial h}=1$ and $\frac{dx}{dt}=H^1_y$. Hence,
 $$\frac{d}{dh}\Big(\int_{\widehat{A_1A_2}}\tilde{q}^1(x,y,\delta)dx\Big)=\int_{\widehat{A_1A_2}}(p^1_x+q^1_y)dt.$$
Further,
$$\begin{aligned}
&\frac{d}{dh}\Big(\int_{\widehat{A_2O}}{q}^1(x,0,\delta)dx\Big)=\frac{d}{dh}\Big(\int_{a_2(h)}^0{q}^1(x,0,\delta)dx\Big)
=-q^1(a_2(h),0,\delta)a'_2(h),\\
&\frac{d}{dh}\Big(\int_{\widehat{O A_1}}p^1(0,y,\delta)dy\Big)=\frac{d}{dh}\Big(\int^{a_1(h)}_0p^1(0,y,\delta)dy\Big)
=p^1(0,a_1(h),\delta)a'_1(h),\\
\end{aligned}$$
which implies the desired conclusion. This ends the proof.\quad$\lozenge$
\end{document}